\documentclass [11pt]{article}
\usepackage{amsthm}
\usepackage{amssymb}
\usepackage{epsfig}
\usepackage{amsthm}
\usepackage{srcltx}
\usepackage[dvipsnames,usenames]{color}
\usepackage{paralist}
\usepackage{amsmath}
\usepackage{amsfonts}

\newcounter{contador}
\newtheorem{propo}[contador]{Proposition}
\newtheorem{teo}[contador]{Theorem}

\newtheorem{remark}[contador]{Remark}
\newcommand{\pr}{{P}{\mathbb C^2}}
\newcommand{\pru}{{P}{\mathbb C^1}}

\newcommand{\R}{{\mathbb R}}
\newcommand{\C}{{\mathbb C}}
\newcommand{\N}{{\mathbb N}}

\newcommand{\Z}{{\mathbb Z}}
\newcommand{\al}{\alpha}
\newcommand{\be}{\beta}
\newcommand{\ga}{\gamma}
\newcommand{\om}{\omega}

\newcommand{\homog}{[x_0:x_1:x_2]}
\title{Invariant Fibrations for some Birational Maps of ${\mathbb C}^2$\footnote{{\bf Acknowledgements.}
 The
authors are supported by Ministry of Economy and Competitiveness of
the Spanish Government through grants MTM2013-40998-P. They are also
supported by the grant  2014-SGR-568 from AGAUR, Generalitat de
Catalunya.
 }}

\author{Anna Cima$^{(1)}$  and Sundus Zafar$^{(1)}$
  \\*[.1truecm]
{\small \textsl{$^{(1)}$ Departament de Matem\`{a}tiques, Facultat
de Ci\`{e}ncies,}}
\\*[-.25truecm] {\small \textsl{Universitat Aut\`{o}noma de Barcelona,}}
\\*[-.25truecm] {\small \textsl{08193 Bellaterra, Barcelona, Spain}}
\\*[-.25truecm] {\small \textsl{cima@mat.uab.cat,
sundus@mat.uab.cat}}}
\date{}
\begin{document}

% ********************** EN CAS D'ARTICLE *********************
\maketitle
\begin{abstract}
%In this article we extract the zero entropy mappings included in a
%certain family of birational mappings of the plane and we give the
%invariant fibrations associated to them.
In this article we extract and study the zero entropy subfamilies of
a certain family of birational maps of the plane. We find these zero
entropy mappings and give the invariant fibrations associated to
them.
\end{abstract}
% *************************************************************

\noindent {\sl  Mathematics Subject Classification 2010:} 14E05,
26C15, 34K19, 37B40, 37C15, 39A23.

\noindent {\sl Keywords:} Birational maps, Algebraic entropy, First
Integrals, Fibrations, Blowing-up, Integrability, Periodicity.

\section{Introduction}\label{S-intro}
A mapping $f=(f_1,f_2):\C^2\longrightarrow \C^2$ is said to be
rational if each coordinate function is rational, that is, $f_i$ is
a quotient of polynomials for $i=1,2.$ These maps can be naturally
extended to the projective plane $\pr$ by considering the embedding
$(x_1,x_2)\in\C\to [1:x_1:x_2]\in\pr.$ The induced mapping
$F:\pr\longrightarrow \pr$ has three components $F_i[x_0:x_1:x_2]$
which are homogeneous polynomials of the same degree. If
$F_1,F_2,F_3$ have no common factors and have degree $d,$ we say
that $f$ or $F$ has degree $d.$ Similarly we can define the degree
of ${F^n} = F \circ \cdot \cdot \cdot  \circ F$ for each $n\in\N.$

We are interested in birational maps. It is said that a rational
mapping $f:\C^2\longrightarrow \C^2$ is birational  if it exists an
algebraic curve and another rational map $g$ such that $f\circ g=g
\circ f =id$ in $\C^2\setminus V.$

The study of the dynamics generated by birational mappings in the
plane has been growing in recent years, see for instance
\cite{BK2,BGM,CGM1,CGM3,JD,IR,JR,LG,LGK,LKP,LKMR,ZE}.

It can be seen that if $f(x_1,x_2)$ is a birational map, then the
sequence of the degrees of ${F^n}$ satisfies a homogeneous linear
recurrence with constant coefficients (see \cite{DF} for instance).
This is governed by the characteristic polynomial $\mathcal{X}(x)$
of a certain matrix associated to $F.$ The other information we get
from $\mathcal{X}(x)$ is the \textit{dynamical degree}, $\delta(F),$
which is defined as
\begin{equation}\label{eq.1}
\delta(F): = \mathop {\lim }\limits_{n \to \infty } {\left( {\deg
({F^n})} \right)^{\frac{1}{n}}}.
\end{equation}
The logarithm of this quantity has been called the \textit{algebraic
entropy of $F$}. It is known that the algebraic entropy is an upper
bound of the topological entropy, which in turn is a dynamic measure
of the complexity of the mapping. For instance, periodic or
integrable birational mappings have zero algebraic entropy.

Birational mappings with zero algebraic entropy have been
characterized, see \cite{DF} and \cite{BD}. From its results we know
the existence of some fibrations associated to the mapping, which
give almost a complete dynamical information of the mapping.

In this paper we consider the family of fractional maps $f:\C^2 \to
\C^2$ :
\begin{equation}\label{ec1}
f(x,y) = \left( {\alpha _0} + {\alpha _1}x + {\alpha
    _2}y,\frac{{\beta _0} + {\beta _1}x + {\beta _2}y}{{\gamma _0} +
    {\gamma _2}y} \right)\,,\,\alpha_1\ne
0\,,\,\beta_1\ne 0\,,\,\gamma_2\ne 0.
\end{equation}
This family is part of a more general family studied in \cite{CZ1}
and \cite{CZ2}, which in turn is a generalization of the birational
mappings studied by Bedford and Kim in \cite{BK2}. The goal of this
paper is to extract, under affine equivalence, all mappings of type
(\ref{ec1}) having zero algebraic entropy and give the corresponding
invariant fibrations associated to them.

The methodology involves the implementation of the blowing-up
technique and the extension of the mappings at the Picard group (see
Section 2).

In general, given a parametric family of mappings, to decide for
which values of the parameters the mappings are periodic, is not an
easy problem (see \cite{CL}, \cite{CGM}, for instance). When the
mapping is a plane birational mapping it is possible to face that
problem (see \cite{BK2}) and it is fascinating to see how these
cases arise, and not only the periodic ones, also all the zero
entropy cases.

The paper is organized as follows. In Section 2 we give some
preliminary results and we explain how we proceed to find the
invariant fibrations associated to zero entropy maps. Section 3
deals with the subfamily $\al_2\ne 0.$ The main result is Theorem
\ref{section1}. Similarly in Section 4 we consider the subfamily
$\al_2=0$ getting Theorem \ref{section2}.

\section{Preliminary results}

Rational mappings $F:{\pr} \to {\pr}$ have an indeterminacy set
${\mathcal{I}(F)}$ of points where $F$ is ill-defined as a
continuous map. This set is given by:
\begin{footnotesize}
$${\mathcal{I}(F)}=\{\homog\in\pr:F_1\homog=0\,,\,F_2\homog=0\,,\,F_3\homog=0]\}.$$
\end{footnotesize}
If $F$ is birational then we can also consider the indeterminacy
points of its inverse $F^{-1}.$ On the other hand, if we consider
one irreducible component $V$ of the determinant of the Jacobian of
$F$, it is known (see Proposition $3.3$ in \cite{JD}) that $F(V)$
reduces to a point which belongs to $\mathcal{I}(F^{-1})$. The set
of these curves which are sent to a single point is called the {\it
exceptional locus} of $F$ and it is denoted by $\mathcal{E}(F).$

It is known that the dynamical degree depends on the orbits of the
indeterminacy points of the inverse of $F$ under the action of $F,$
see \cite{DF, FS}. Indeed, the key point is whether the iterates of
such points coincide with any of the indeterminacy points of $F.$
When it happens, this orbit is finite.

Sometimes some \textit{orbit collision} appears. The expression
orbit collision refers to the following: Let $S\in\mathcal{E}(F)$
which collapses at the point $A\in \mathcal{I}(F^{-1})$ (we will
write $S\twoheadrightarrow A$ to describe this behaviour). Following
the orbit of $A,$ assume that it ends at a point $O\in
\mathcal{I}(F).$ It can happen that
$$S\twoheadrightarrow A\rightarrow \ast \rightarrow\cdots
\rightarrow \sigma  \rightarrow\cdots \rightarrow O \,,\,\sigma\in
\bar{S}\in\mathcal{E}(F).$$ Then, being $f$ birational, (see
\cite{DF}) it exists $\bar{A}\in\mathcal{I}(F^{-1})$ with
$\bar{S}\twoheadrightarrow \bar{A}.$ When it happens it is said that
the orbits of $A$ and $\bar{A}$ collides. This is exactly the
behaviour that we get in family (\ref{eq.1}) and what makes the
family so interesting.

\subsection{Blow-up's and the Picard group}
Given a point $p\in\C^2,$ let $(X,\pi),$ be the blowing-up of $\C^2$
at the point $p.$ Then,translating $p$ at the origin,
$$\pi^{-1}p=\pi^{-1}(0,0)=\{\left((0,0),[u:v]\right)\}:=E_p\simeq\pru$$
and if $q=(x,y)\ne (0,0),$ then
$$\pi^{-1}q=\pi^{-1}(x,y)=\left((x,y),[x:y]\right)\in X.$$
Given the point $\left((0,0),[u:v]\right)\in E_{p}$ (resp.
$\left((x,y),[x:y]\right)$) we are going to represent it by
$[u:v]_{E_p}$ (resp. by $(x,y)\in\C^2$ or by $[1:x:y]\in\pr$ if it
is convenient). After every blow up we get a new expanded space $X$
and the induced map $\tilde{F}: X \to X.$ And then $\tilde F$
induces a  morphism of groups,
 $\tilde F^*:\mathcal{P}ic(X)\rightarrow \mathcal{P}ic(X)$ just by
 taking classes of preimages, where $\mathcal{P}ic(X)$ is the Picard group of $X$ (see \cite{BK1,BK2}).
 It is proved that after a finite number of blowing-up's
 we get a map $\tilde{F}$ which satisfies
 $\left(\tilde{F}^n\right)^*=\left(\tilde{F}^*\right)^n.$ Maps $\tilde F$ satisfying
 this equality are called \textit{Algebraically Stable Maps} (AS for
 short), (see \cite{DF}). The characteristic polynomial of the
 matrix of $\tilde{F}^*$ is the one associated to the sequence of
 degrees $d_n:=\,\text{degree}\,F^n.$

\subsection{Lists of orbits.}
We derive our results by using Theorem \ref{th_BK} below,
established and proved in \cite{BK1,BK2}. The proof of that is based
in the same tools explained in the above paragraph. In order to
determine the matrix of the extended map in the Picard group, it is
necessary to distinguish between different behaviors of the iterates
of the map on the indeterminacy points of its inverse.

The theorem is written for a general family $G$ of quadratic maps of
the form $G = L \circ J.$ As we will see the maps of  family
(\ref{eq1}), when the triangle is non-degenerate, are linearly
conjugated to such a maps. Here $L$ is an invertible linear map and
$J$ is the involution in $\pr$ as follows:
\begin{equation*}
J[x_0:x_1:x_2] = [{x_1}{x_2}:{x_0}{x_2}:{x_0}{x_1}].
\end{equation*}
We find that the involution $J$ has an indeterminacy locus
$\mathcal{I} = \{\epsilon_0, \epsilon_1, \epsilon_2\}$ and a set of
exceptional curves  $\mathcal{E} = \{\Sigma_0, \Sigma_1,
\Sigma_2\}$, where $\Sigma_i = \{x_i = 0\}$ for $i = 0,1,2,$ and
$\epsilon_i = \Sigma_j \bigcap \Sigma_k$ with $\{i,j,k\} =\{0,1,2\}$
and $i\ne j\ne k\,,\,i\ne k.$ Let $\mathcal{I}(G^{-1}) := \{a_0,
a_1, a_2\},$ the elements of this set are determined by $a_i:=
G(\Sigma_i - \mathcal{I}(J))=L \,\epsilon_i$ for $i = 0,1,2;$ see
\cite{BK1}.

To follow the orbits of the points of $\mathcal{I}(G^{-1})$ we need
to understand the following definitions and construction of lists of
orbits in order to apply the result of Theorem \ref{th_BK}.

We assemble the orbit of a point $p \in \pr$ under the map $G$ as
follows. For a point $p \in \mathcal{E}(G) \cup \mathcal{I}(G)$  we
say that the orbit $\mathcal{O}(p) = \{p\}$. Now consider that there
exits a $p \in \pr$ such that its $n^{th}-$ iterate belongs to
$\mathcal{E}(G) \cup \mathcal{I}(G)$ for some $n$, whereas all the
other $n-1$ iterates of $p$ under $G$ are never in $\mathcal{E}(G)
\cup \mathcal{I}(G)$. This is to say that for some $n$ the orbit of
$p$ reaches an exceptional curve of $G$ or an indeterminacy point of
$G.$ We thus define the orbit of $p$ as $\mathcal{O}(p) = \{p, G(p),
..., G^{n}(p)\}$ and we call it a  \textit{singular orbit}. If for
some $p \in \pr$ in turns out that $p$ and all of its iterates under
$G$ are never in $\mathcal{E}(G) \cup \mathcal{I}(G)$ for all $n,$
we set as $\mathcal{O}(p) = \{p, G(p), G^{2}(p) ...\}$ and
$\mathcal{O}(p)$ is \textit{non singular orbit}. We now make another
characterization of these orbits. Consider that a singular orbit
reaches an indeterminacy point of $G$, this is to say that $G^{n}(p)
\in \mathcal{I}(G)$ but its not in $\mathcal{E}(G).$  We call such
orbits as \textit{singular elementary orbits} and we  refer them as
SE-orbits. To apply Theorem \ref{th_BK} we need to organize our SE
orbits into lists in the following way.

Two orbits $\mathcal{O}_{1} = \{a_1,...,\epsilon_{j_1}\}$ and
$\mathcal{O}_{2} = \{a_2,...,\epsilon_{j_2}\}$ are in the same list
if either $j_1=2$ or $j_2=1,$ that is, if the ending index of one
orbit is the same as the beginning index of the other. We say that a
list of orbits
$$\mathcal{L} = \left\{\mathcal{O}_{i} = \{a_i,...,\epsilon_{\tau(i)}\},\ldots,\mathcal{O}_{j} = \{a_j,...,\epsilon_{\tau(j)}\}\right\}$$
is closed if $\tau(j)=i.$ Otherwise it is an open list. For
instance,
$$\begin{array}{lll}
\mathcal{L}_1 = \left\{\mathcal{O}_{1} = \{a_1,...,\epsilon_1\}\right\},\\
\mathcal{L}_2 = \left\{\mathcal{O}_{0} = \{a_0,...,\epsilon_{2}\},\mathcal{O}_{2} = \{a_2,...,\epsilon_{0}\}\right\},\\
\mathcal{L}_3 = \{\mathcal{O}_{0} = \{a_0,...,\epsilon_1\},
\mathcal{O}_{1} = \{a_1,...,\epsilon_2\}, \mathcal{O}_{2} =
\{a_2,...,\epsilon_0\}\},
\end{array}$$
are closed lists.

We now define two polynomials $\mathcal{T}_{\mathcal{L}}$ and
$\mathcal{S}_{\mathcal{L}}$ which we will use to state Theorem
\ref{th_BK}. Let $n_i$ denote the sum of the number of elements of
an orbit $\mathcal{O}_{i}$ and let $\mathcal{N}_{\mathcal{L}} = n_u
+ ...+ n_{u+\mu}$ denote the sum of the numbers of elements of each
list $\left|\mathcal{L}\right|.$ If $\mathcal{L}$ is closed then
$\mathcal{T}_{\mathcal{L}} = x^{\mathcal{N}_{\mathcal{L}}}-1$ and if
$\mathcal{L}$ is open then $\mathcal{T}_{\mathcal{L}} =
x^{\mathcal{N}_{\mathcal{L}}}.$ Now we define
$\mathcal{S}_{\mathcal{L}}$ for different lists as follows:
\begin{equation*}
{\mathcal{S}_{\mathcal{L}}}(x) = \left\{ {\begin{array}{*{20}{c}}

{1} & {} & \mbox{if}\;{\left| \mathcal{L} \right| = \left\{ {{n_1}} \right\},}  \\
   {{x^{{n_1}}} + {x^{{n_2}}} + 2} & {} & \mbox{if}\; \mathcal{L}\; \mbox{is closed and}\; {\left| \mathcal{L} \right| = \left\{ {{n_1},{n_2}} \right\},}  \\
   {{x^{{n_1}}} + {x^{{n_2}}} + 1} & {} & \mbox{if}\;\mathcal{L}\; \mbox{is open and}\; {\left| \mathcal{L} \right| = \left\{ {{n_1},{n_2}} \right\},}  \\
   {\sum\limits_{i = 1}^3 {\left[ {{x^{{\mathcal{N}_{\mathcal{L}}} - {n_i}}} + {x^{{n_i}}}} \right]  + 3}} & {} & \mbox{if}\; \mathcal{L}\; \mbox{is closed and}\; {\left| \mathcal{L} \right| = \left\{ {{n_1},{n_2},{n_3}} \right\},}  \\
   {\sum\limits_{i = 1}^3 {{x^{{\mathcal{N}_{\mathcal{L}}} - {n_i}}}}  + \sum\limits_{i \ne 2} {{x^{{n_i}}}}  + 1} & {} & \mbox{if}\;\mathcal{L}\; \mbox{is open and}\;{\left|\mathcal{L} \right| = \left\{ {{n_1},{n_2},{n_3}} \right\}.}  \\
 \end{array} } \right.
\end{equation*}
\begin{teo}\label{th_BK}
(\cite{BK2})\,\,If $G = L \circ J$, then the dynamical degree
$\delta(G)$ is the largest real zero of the polynomial
$$\mathcal{X}(x) = (x - 2)\prod\limits_{\mathcal{L} \in {\mathcal{L}^c} \cup {\mathcal{L}^o}} {{\mathcal{T}_{\mathcal{L}}}(x) + (x - 1)\sum\limits_{\mathcal{L} \in {\mathcal{L}^c} \cup {\mathcal{L}^o}} {{S_L}(x)\prod\limits_{\mathcal{L}' \ne \mathcal{L}} {{\mathcal{T}_{\mathcal{L}'}}(x).} } }$$
Here $\mathcal{L}$ runs over all the orbit lists.
\end{teo}

This theorem enables us to calculate the characteristic polynomial
associated to $d_n.$ To this end we have to perform the lists of the
orbits of the points in $\mathcal{I}(F^{-1}),$ but for this we have
to do the necessary blow-up's to get an AS mapping.

In order to get AS maps we will use the following useful result
showed by Fornaess and Sibony in \cite{FS} (see also Theorem 1.14)
of \cite{DF}:
\begin{equation}\label{condicio}
{\text{The map $\tilde{F}$ is AS if and only if for every
exceptional curve}\,\, C\,\, \text{and all} \,\, n\ge 0\,,\,\tilde
F^n(C)\notin {\mathcal{I}}(\tilde F).}
\end{equation}

\subsection{Zero entropy}
The following result is quiet useful in our work. It is a direct
consequence of Theorem $0.2$ of \cite{DF}. Given a birational map
$F$ of $\pr,$ let $\tilde{F}$ be its regularized map so that the
induced map $\tilde{F}^{*}:\mathcal{P}ic(X) \to \mathcal{P}ic(X)$
satisfies $(\tilde{F}^n)^{*} = (\tilde{F}^{*})^n.$ Then

\begin{teo}\label{theo-diller}
    (See \cite{DF}) Let $F:\pr \to \pr$ be a birational map, $\tilde{F}$ be its regularized map and let $d_n = deg(F^n).$  Then up to bimeromorphic conjugacy, exactly one of
    the following holds:
    \begin{itemize}
        \item The sequence $d_n$ grows quadratically, $\tilde{F}$ is an automorphism and $f$ preserves an elliptic fibration.
        \item The sequence $d_n$ grows linearly and $f$ preserves a rational fibration. In this case $\tilde{F}$ cannot be conjugated to an automorphism.
        \item The sequence $d_n$ is bounded, $\tilde{F}$ is an automorphism and $f$ preserves two generically transverse rational fibrations.
        \item The sequence $d_n$ grows exponentially.
    \end{itemize}
    In the first three cases $\delta(F) = 1$ while in the last one $\delta(F) > 1.$ Furthermore in the first and second, the invariant fibrations are unique.
\end{teo}

We recall that $f:\C^2 \to \C$ \textit{preserves a fibration}
$V:\C^2\to \C$ if $f$ sends level curves of $V$ to level curves of
$V.$ If $f$ sends each level curve of $V$ to itself, it is said that
$f$ is \textit{integrable} and that $V$ is a first integral of $f.$

When the sequence $d_n$ is bounded it can happen that it is periodic
or not. For mappings which are not periodic, we have the following
result (Theorem A of \cite{BD}):

\begin{teo}\label{blancdeserti}(See \cite{BD})
Let $F:\pr\rightarrow \pr$ be a non-periodic birational map such
that the corresponding sequence of degrees is bounded. Then $F$ is
conjugate to an automorphism of $\pr,$ which restricts to one of the
following automorphisms on some open subset isomorphic to $\C^2:$
\begin{itemize}
\item [(1)] $(x,y)\mapsto (\al\,x,\beta\,y),$ where
$\al,\beta\in\C^*,$ and where the kernel of the group homomorphism
$\Z^2\to \C^*$ given by $(i,j)\mapsto \al�\,\beta^j$ is generated by
$(k,0)$ for some $k\in\Z.$
\item [(2)] $(x,y)\mapsto (\al\,x,y+1),$ where $\al\in\C^*.$
\end{itemize}
\end{teo}

\subsection{Invariant fibrations}\label{invfib}
From Theorem \ref{theo-diller} we know the existence of rational
invariant fibrations depending on the growth of $d_n.$ To find them,
we consider $V(x,y)=\frac{P(x,y)}{Q(x,y)}$ for some polynomials
$P(x,y),Q(x,y)$ without common factors. If $V$ is an invariant
fibration, then $f$ sends $V=k$ to $V=k'.$ In this work we consider
that the relation between $k,k'$ is of type
$\psi(k)=\frac{\om_1\,k+\om_2}{\om_3\,k+\om_4}$ for some
$\om_1,\om_2,\om_3,\om_4\in\C.$ In particular we will have the
following cases:
\begin{itemize}
\item [(a)] $V(f)=V,$ the integrable case.
\item [(b)] $V(f)=\om_1\,V,$ the scaled fibration case.
\item [(c)] $V(f)=\om_1\,V+\om_2,$ the scaled translated fibration
case. \end{itemize} Note that in case $(a)$ in general the functions
$P$ and $Q$ are invariant under $f$ as they satisfy the equation
$P\cdot Q(f)=Q\cdot P(f)$ (unless that the denominators of $P(f)$ or
$Q(f)$ are simplified with $Q$ or $P$ respectively). Similarly for
case $(b)$ it follows. In case $(c)$ only $Q$ is invariant as it
satisfies the relation $Q\cdot P(f)=(\om_1\,P+\om_2\,Q)\cdot Q(f).$

Hence we always begin finding invariant algebraic curves. To find
them, we introduce the following definition. Given a birational map
and given a curve $C\subset \pr$ we define
$F(C):=\overline{F(C\setminus \mathcal{I}(F)}$ to be the proper
transform of $C$ by $F.$ When $C\cap \mathcal{I}(F)=\emptyset,$ we
have that deg $F(C)=d\cdot$ deg$(C)$ where $d$ is the degree of $F.$
In general,
\begin{equation}\label{inv}
\text{deg}\,F(C)\,=\,d\cdot\text{deg}\,(C)-\sum_{O\in\mathcal{I}(F)}m_O(C),
\end{equation}
where $m_O(C)$ is the algebraic multiplicity of $C$ at $O$ (see (1),
pg. 416, \cite{D2}).

The approach is the following. Take an arbitrary curve $C$ and
impose that deg$\,F(C)= $ deg$\,C,$ that is,
$(d-1)\,\text{deg}\,(C)=\sum_{O\in\mathcal{I}(F)}m_O(C).$ For
instance if $d=2$ and we consider $C$ of degree $3,$ then a
necessary condition for $C$ to be invariant under $f$ is that $C$
passes through three indeterminacy points of $F$ of multiplicity one
or through one indeterminacy point with multiplicity two and another
of multiplicity one or through one indeterminacy point with
multiplicity three. In the first case, for instance, if
$O_1,O_2,O_3\in\mathcal{I}(F),$ then there exist
$T_1,T_2,T_3\in\mathcal{E}(F)$ such that
$F^{-1}:T_i\twoheadrightarrow O_i.$ Also, if $C=\{P=0\}$ for some
polynomial $P,$ then $F(C)\subset\{P\circ F^{-1}=0\}$ and we have
that $P\circ F^{-1}=T_1\cdot T_2\cdot T_3\cdot \bar{P}$ for a
certain polynomial $\bar{P},$ with $F(C)=\{\bar{P}=0\}.$ Then
imposing that $P-k\cdot \bar{P}=0$ we will find, if we get a
solution, an invariant curve of degree three.

As we will see our particular mappings, sometimes depend on a number
$\al$ which is a zero of certain polynomial $P$. Then all the
calculations have to be made in $\frac{\C[\al]}{(P(\al))}[x,y],$
which in fact make them more complicated ($\frac{\C[\al]}{(P(\al))}$
is the quotient ring $\C[\al]$ over the ideal generated by the
polynomial $P(\al)$).

\section{The subfamily  $\alpha_2 \ne 0.$}

Taking into account that $\alpha_1,\beta_1$ and $\gamma_2$ are not
zero, it can be proved  that when $\alpha_2 \ne 0,$ after an affine
change of coordinates $f(x,y)$ can be written as
\begin{equation}\label{eq2}f(x,y)=\left(\alpha_0+\alpha_1 x+y,\frac{x}{\gamma_0+y}\right)\quad ,\quad \alpha_1\ne
0.\end{equation}

We consider the imbedding $(x,y) \mapsto [1:x:y] \in {\pr}$ into
projective space and consider the induced map $F:{\pr} \to {\pr}$
given by:
\begin{equation}\label{F3ND:Eq2}
F[{x_0}:{x_1}:{x_2}] = [x_0(\ga_0 x_0 + x_2):(\al_0 x_0 + \al_1 x_1
+ x_2)(\ga_0 x_0 + x_2): x_0 x_1].
\end{equation}
The indeterminacy locus of $F$ is $\mathcal{I}(F) = \left\{
{{O_0},{O_1},{O_2}} \right\},$ where
\begin{equation*}
{O_0} = \left[ 1:0:{-\ga_0}\right],\quad {O_1} = \left[ { 0: 1:
-\al_1} \right],
 \quad {O_2} = \left[ {0: 1: 0} \right],
\end{equation*}
and the indeterminacy locus of $F^{-1}$ is $\mathcal{I}(F^{-1}) =
\left\{ {{A_0},{A_1},{A_2}} \right\},$ where
\begin{equation*}
{A_0} = \left[ 0 : 1 : 0\right],\quad {A_1} = \left[ { 0: 0 : 1}
\right], \quad {A_2} = \left[ -\al_1: -\al_1 (\al_0- \ga_0): 1
\right].
\end{equation*}
The set of exceptional curves is given as $\mathcal{E}(F) = \left\{
{{S_0},{S_1},{S_2}} \right\},$ where
\begin{equation*}
{S_0} = \left\{ {{x_0} = 0} \right\},\quad {S_1} = \left\{ {\ga_0
x_0 + x_2 = 0} \right\},\quad {S_2} = \left\{ {\ga_0 x_0 + x_2 +
\al_1 x_1 = 0}\right\},
\end{equation*}
and the set of exceptional curves of $F^{-1}$ is given as
$\mathcal{E}(F^{-1}) = \left\{ {{T_0},{T_1},{T_2}} \right\},$ where
\begin{equation*}
{T_0} = \left\{ {(\al_0 - \ga_0)x_0 - x_1 = 0} \right\},\quad {T_1}
= \left\{ x_0 - \al_1 x_2  \right\},\quad {T_2} = \left\{ {x_0 =
0}\right\}.
\end{equation*}

\begin{teo}\label{theo2}
Let $F[x_0:x_1:x_2]$ be defined by
$$F[x_0:x_1:x_2]=[x_0(\ga_0x_0+x_2):(\al_0x_0+\al_1x_1+x_2)(\ga_0x_0+x_2):x_0x_1]$$
and let $\tilde{F}$ be the induced map after blowing up the point
$A_0.$ Then the following hold:
\begin{itemize}
    \item If $\tilde{F}^{2k}(A_1) \ne
O_1$ for all $k \in  \N$ and $\tilde{F}^{p}(A_2) = O_0$ for some $p \in \N$ then the characteristic polynomial
associated with $F$ is given by
\begin{equation*}
\mathcal{X}_p = x^{p+1}(x^2-x-1)+x^2,
 \end{equation*}
 and
 \begin{itemize}
    \item for $p = 0,\,p=1$ the sequence of degrees $d_{n}$ is bounded,
    \item for $p = 2$ the sequence of degrees $d_{n}$ grows linearly,
    \item for $p > 2$ the sequence of degrees $d_{n}$ grows exponentially.
\end{itemize}

 \item Assume that $\tilde{F}^{2k}(A_1) = O_1$ for some $k \in \N.$ Let
$\tilde{F}_1$ be the induced map after we blow-up the points
$A_0,A_1,\tilde{F}(A_1),\ldots ,\tilde{F}^{2k}(A_1)=O_1.$  If
$\tilde{F}_1^{p}(A_2) \ne O_0$ for all $p\in\N,$ then the
characteristic polynomial associated with $F$ is given by

\begin{equation*}\label{sik}
\mathcal{X}_k = x^{2k+1}(x^2-x-1)+1,
\end{equation*}
and the sequence of degrees grows exponentially. Furthermore
$\delta(F) \to \delta^{*}$ as $k \to \infty.$

\item If $\tilde{F}^{2k}(A_1) = O_1$ and $\tilde{F}^{p}(A_2) = O_0$ for some $p,k \in \N$ then the characteristic polynomial
associated with $F$ is given by
\begin{equation*}
\mathcal{X}_{(k,p)} =
x^{p+1}(x^{2k+3}-x^{2k+2}-x^{2k+1}+1)+x^{2k+3}-x^2-x+1,
 \end{equation*}
and
\begin{itemize}
    \item for $p > \frac{2\,(1+k)}{k}$ the sequence of degrees $d_{n}$ grows exponentially.
    \item for $(k,\,p) \in \{(2,\,3),\,(1,\,4)\}$ the sequence of degrees $d_{n}$ either, it is periodic or it grows quadratically;
    \item for $(k,\,p) \in \{(k,\,0),\,(k,\,1),\,(k,\,2),\,(1,\,3)\}$ the sequence of degrees $d_{n}$ is periodic.
    \end{itemize}
 \item Assume that $\tilde{F}^{2k}(A_1)\ne O_1$ and $\tilde{F}^p(A_2)\ne O_0$ for all
$k,\,p\in\N.$ Then the characteristic polynomial associated with $F$
is given by $$\mathcal{X}(x) = {x^2} - x - 1,$$ and the sequence of
degrees grows exponentially with $\delta(F)=\delta^{*}.$
\end{itemize}
\end{teo}

\begin{proof}
Observe that $S_0 \twoheadrightarrow A_0 = O_2.$ The orbit of $A_0$
is SE. By blowing up $A_0$ we get the exceptional fibre $E_0$ and
the new space $X.$ The induced map $\tilde{F}: X \to X$ sends the
curve $S_0 \to E_0 \to S_0.$ Observe that now
$\mathcal{{I}}(\tilde{F}) = \left\{ {{O_0},{O_1}} \right\}$ and
$\mathcal{{E}}(\tilde{F}) = \left\{{{S_1},{S_2}} \right\}.$

We see that $A_1 \ne O_1$ and the exceptional curve $S_1
\twoheadrightarrow A_1 \in S_0.$ We observe that the collision of
orbits discussed in preliminaries is happening here. The orbit of
$A_1$ under $\tilde{F}$ is as follows:
\begin{equation*}
S_1 \twoheadrightarrow  A_1 \to [\gamma_2 :\beta_2]_{E_0} \to [0 :
\alpha_1(\gamma_0+\beta_2): \beta_1]\in S_0 \to \cdots
\end{equation*}
After some iterates we can write the expression of
$\tilde{F}^{2k}(A_1)$ for all $k > 0 \in \N$ as
%\begin{equation*}\label{eq3}
$\tilde{F}^{2k}(A_1) = [0 :
\alpha_1(\gamma_0+\beta_2)(1+\alpha_1+\alpha_1^2+\cdots+\alpha_1^{k-1}):
\beta_1]\in S_0.$
%\end{equation*}
Observe that for some value of $k \in \N$ it is possible that
$\tilde{F}^{2k}(A_1) = O_1.$ This happens when the following
\textit{condition $k$}
 is satisfied for some $k.$
\begin{equation}\label{eq4}
 \alpha_1^2(\gamma_0+\beta_2)(1+\alpha_1+\alpha_1^2+\cdots+\alpha_1^{k-1})+ \alpha_2\beta_1=0.
\end{equation}
For such $k \in \N$ the orbit of $A_1$ is SE. By blowing up the
points of this orbit we get the new space $X_1$ and the induced map
$\tilde{F}_1.$ Then under the action of $\tilde{F}_1$ we have
\begin{equation*}\label{seq2} S_1 \to G_0 \to G_1 \to G_2 \to \cdots
\to G_{2k-1}\to G_{2k}\to T_1.\end{equation*} Then
$\mathcal{I}(\tilde{F}_1)=\{O_0\}$ and
$\mathcal{E}(\tilde{F}_1)=\{S_2\}.$

Now if the orbit of $A_1$ is SE and if $\tilde{F}_1^{p}(A_2) = O_0$
that is the orbit of $A_2$ is also SE for some $p \in \N$ then we
have three SE orbits. If \textit{condition $k$} is not satisfied
then with the extended map $\tilde{F}$ we have
$\mathcal{I}(\tilde{F})=\{O_0,O_1\}.$ Therefore we have two options:
$\tilde{F}^{p}(A_2) = O_0$ or $\tilde{F}^{p}(A_2) = O_1.$

We claim that for all $p\in\N,$ $\tilde{F}^{p}(A_2) \ne  O_1.$
Assume that $\tilde{F}^{p}(A_2)= O_1$ and assume that
$F^j(A_2)\notin S_0$ for $j=1,2,\ldots, p-1.$
$\tilde{F}^{p}(A_2)={F}^{p}(A_2)=O_1.$ Since $O_1\in S_0$ and
$A_2\notin S_0$ if $F^p(A_2)= O_1$ then $p$ would be greater than
zero and since $S_0=T_2,$ it would imply that $O_1=A_1$ or
$O_1=A_2,$ which is not the case (recall that the only points in
$T_2$ which have a preimage are $A_1$ and $A_2$).

Contrarily, if it exists some $l\in\N,l<p$ such that $F^j(A_2)\notin
S_0$ for $j=1,2,\ldots, l-1$ but $F^l(A_2)\in S_0\setminus\{O_1\}$
then $F^l(A_2)$ must be equal to $A_1$ or $A_2$ that is,
${F}^{l}(A_2)= A_1$ or ${F}^{l}(A_2)= A_2.$ The second case is not
possible as $A_2$ is a fixed point. In the first case
$\tilde{F}^{p}(A_2)=\tilde{F}^{p-l}(F^l(A_2))=\tilde{F}^{p-l}(A_1)=O_1$
which implies that $p=l+2r$ and $\tilde{F}^{2r}(A_1)=O_1.$  Hence
the orbit of $A_1$ must be SE and that condition $k$ must be
satisfied for $k=r$  which is a contradiction. It implies that the
only available possibility for $\mathcal{O}_2$ to be SE is to have
that for some $p,$ $\tilde{F}^{p}(A_2) = O_0.$ After the blow up
process we get
$$S_2 \to E_1 \to E_2 \to \cdots\to E_{p}\to E_{p+1} \to T_0.$$
The extended map $\tilde{F}_2$ is an automorphism when we have three
SE orbits.

The above discussion gives us three different cases.
\begin{itemize}
    \item One $SE$ orbit: This happens when $A_0 = O_2$ with the conditions that $\tilde{F}^{2k}(A_1)\ne O_1$ and
    $\tilde{F}^p(A_2)\ne O_0$ for all $k,p\in\N.$
    Therefore we have only one list $\mathcal{L}_{o}$ which is open that is
    $\mathcal{L}_{o} = \{ \mathcal{O}_{0} = \{A_0 = O_2\}\}. $ By using Theorem \ref{th_BK} we find that $\delta(F)= \frac{\sqrt{5}+1}{2},$
    which is given by the greatest root of the
    polynomial $X(x) = {x^2} -x- 1.$ Therefore it has exponential growth.
\item Two $SE$ orbits $(a)$: It is the case when $A_0 = O_2, \,\,\tilde{F}^{p}(A_2) = O_0$ and $\tilde{F}^{2k}(A_1) \ne O_1$
for all $k\in\N.$ By organizing the orbits into lists we have one
closed list $\mathcal{L}_c = \{\mathcal{O}_0 = \{A_0 = O_2\},\quad
\mathcal{O}_2 =
\{A_2,\,\tilde{F}(A_2)\,,...,\tilde{F}^{p}(A_2)=O_0\}\}.$ By
utilizing Theorem \ref{th_BK} we find that the characteristic
polynomial associated to $F$ is $\mathcal{X}_p =
x^{p+1}(x^2-x-1)+x^2.$ For $p=0$ and $p=1$ the sequence of degrees
satisfies $d_{n+3} = d_{n} $ and  $d_{n+4} = d_{n+3}$ respectively
which corresponds towards boundedness of $d_n$.

For $p = 2$ we get the polynomial $\mathcal{X}_2 = x^2(x+1)(x-1)^2.$
Looking at the first degrees we get that the sequence of degrees is
$d_n = -1+2\,n.$

For $p > 2,$ we observe that
$\mathcal{X}_p(1)=0,\,\mathcal{X'}_p(1)=2-p<0$ and $\lim_{x\to
+\infty} \mathcal{X}_p(x)=+\infty.$ Hence  $\mathcal{X}_p$ always
has a root $\lambda>1$ and the result follows.

\item Two $SE$ orbits $(b)$: When we have $A_0 = O_2, \,\,\tilde{F_1}^{2k}(A_1)
=O_1$ and $\tilde{F_1}^{p}(A_2)\ne O_0$ for all $p\in\N$ then there
is one open and one closed list and $\mathcal{X}_k =
x^{2k+1}(x^2-x-1)+1.$ We observe that for all the values of $k \in
\N\,,\,k\ge 1$ the polynomial $\mathcal{X}_k$ has always a root
$\lambda > 1.$ Therefore $f$ has exponential growth.

\item Three $SE$ orbits: In this case we have $A_0 = O_2,\,\,\tilde{F}^{2k}(A_1) = O_1,\,\,\tilde{F}^{p}(A_2) = O_0,$ for
a certain $p,k\in\N.$ We have two closed lists as follows:
$$\mathcal{L}_c = \{\mathcal{O}_0 = \{A_0 = O_2\},\quad \mathcal{O}_2 = \{A_2,\,\tilde{F}(A_2)\,,...,\tilde{F}^{p}(A_2)=O_0\}\},$$
$$\mathcal{L}_c= \{\mathcal{O}_1 = \{A_1,\,\tilde{F}(A_1)_{E_0}\,,...,\tilde{F}^{2k}(A_1)_{S_0}=O_1\}\}.$$
From Theorem \ref{th_BK} we can write $\mathcal{X}_{(k,p)} =
x^{p+1}(x^{2k+3}-x^{2k+2}-x^{2k+1}+1)+x^{2k+3}-x^2-x+1.$ The map
$\tilde{F}_2$ is an automorphism for all the values $(k,p).$
According to Diller and Favre in \cite{DF} the degree growth of
iterates of an automorphism could be bounded, quadratic or
exponential but it cannot be linear as in such a case the map is
never an automorphism. For this we observe the behavior of
$\mathcal{X}_{(k,p)}$ around $x=1.$ We consider it's Taylor
expansion near $x = 1:$
$$ \mathcal{X}_{(k,p)}(x) = 2(2-kp+2k)(x-1)^2 +
O(\left|x-1\right|^3).$$ Thus $\mathcal{X}_{(k,p)}$ vanishes at $x =
1$ and has a maximum on it if $p>\frac{2(1+k)}{k}.$ Since
$\lim_{x\to +\infty \mathcal{X}_{(k,p)}(x)}=+\infty,$ always exists
a root greater than one. If $p\le \frac{2(1+k)}{k}\,,\,k\ge 1$ then
the pairs $(k,p)$ are in the set: $A_{(k,p)} = \{((k \geq
1),0),\,((k \geq 1),1),\,((k \geq 1),2),\,(1,3),\,(2,3),\,(1,4)\}.$

For $(k,p) = (k,0),$
$\mathcal{X}_{(k,0)}(x)=(x^{2k+2}-1)\,(x-1)\,(x+1),$ and hence the
sequence of degrees is
$$d_{n} = c_0+c_1 n+c_2\,(-1)^n + c_3\,(-1)^n\,n + c_4\,\lambda_1^n
+ c_5\, \lambda_2^n +...+c_{2k+3}\,\lambda_{2k}^n,$$ where $c_i$ are
constants and $\lambda$'s are the roots of polynomial $x^{2k+2}=1$
different from $\pm 1.$ Since $\tilde{F}_2$ is an automorphism for
all ${(k,p)},$ using \cite{DF} we have $c_1 = 0 = c_3.$ This implies
that $d_{2k+2+n} = d_{n},$ i. e., the sequence of degrees is
periodic with period $2k+2.$ The argument for the proof of other
values of $(k,p) \in A_{(k,p)}$ follows accordingly.

%If $k$ is odd then $$d_{n} = l_0+l_1 n+l_2\,(-1)^n + l_3\,\mu_1^n+
%l_4\,\mu_2^n + \cdots+l_{2k+3}\,\mu_{2k+1}^n,$$ where $l_i$ are
%constants and $\mu$'s are the roots of polynomial
%$(x^{k+1}-1)\,(x^k+x^{k-1}+\cdots +x+1).$ Since $\tilde{F}_2$ is an
%automorphism for all ${(k,p)},$ using \cite{DF} we have $c_1 = 0 =
%c_3$ and also $l_1=0.$ This implies that $d_{2k+2+n} = d_{n},$ i.
%e., the sequence of degrees is periodic with period $2k+2.$ The
%argument for the proof of other values of $(k,p) \in A_{(k,p)}$
%follows accordingly.

\end{itemize}
\end{proof}

From the above theorem we see that zero entropy cases only appear
when $\tilde{F}^p(A_2)= O_0,\,\tilde{F}^{2k}(A_1)\neq O_1$ for $p \in \{0,1,2\}$ and $\forall k \in \N$ and when $\tilde{F}^p(A_2)= O_0,\,\tilde{F}^{2k}(A_1) = O_1$ for $(k,\,p) \in \{(k,\,0),\,(k,\,1),\,(k,\,2),\,(1,3),\,(2,3),\,(1,4)\}.$
We are going to study the dynamics of each case separately. Recall that condition $k$ is given by
\begin{equation}\label{cond.k}\al_1^2\,\ga_0\,(1+\al_1+\al_1^2+\cdots
+\al_1^{k-1})+1=0.\end{equation}

The following proposition considers the case when $p=0.$ From the
above theorem we know that if \textit{condition $k$} is not
satisfied the sequence $d_n$ is bounded and when it is satisfied,
$d_n$ is a periodic sequence of period $2k+2.$ In any case we have
to find two generically transverse fibrations. In the second case we
present two first integrals functionally independent. We also prove
that when $d_n$ is periodic, the mapping $f(x,y)$ is itself
periodic.

%and \textit{condition $k$} is not satisfied. From the above theorem
%we know that the sequence $d_n$ is bounded therefore we look for two
%generically transverse fibrations as follows:

\begin{propo}\label{pzero}
Assume that $A_2=O_0.$ Then $f(x,y)$ can be written as
%$f(x,y)$ can Assume that
%$\tilde{F}^{2k}(A_1) \ne O_1$ and $\tilde{F}^{p}(A_2) = O_0$ for
%$(k,\,p) \in \{(k,\,0) : \,k \in \N\},$ then $f(x,y)$ can be written

\begin{equation}\label{pzero}f(x,y)=\left(\frac{1}{\alpha_1}+\alpha_1 x+y,\frac{x}{\frac{1}{\alpha_1}+y}\right)\quad ,\quad \alpha_1\ne
0\end{equation} and the following hold:
\begin{itemize}
\item If $\al_1\ne 1$ then $f(x,y)$ preserves the
two generically transverse fibrations
 $$V_1(x,y)=\frac{\sqrt{\alpha_1}-\alpha_1(\alpha_1+\sqrt{\alpha_1})x+\alpha_1(1+2\sqrt{\alpha_1})y+\alpha_1^2(1+\sqrt{\alpha_1})y^2}{1+\alpha_1
 y}$$
 $$V_2(x,y)=\frac{-1+\alpha_1(1-\sqrt{\alpha_1})x+(\sqrt{\alpha_1}-2\alpha_1)y+\alpha_1(\sqrt{\alpha_1}-\alpha_1)y^2}{1+\alpha_1  y}$$
with $V_1(f(x,y))=-\sqrt{\alpha_1}\,V_1(x,y)$ and
$V_2(f(x,y))=\sqrt{\alpha_1}\,V_2(x,y).$

If $\al_1^{k+1}=1$ then $f$ is a $(2k+2)-$periodic map. In this case
$W_1(x,y)$ and $W_2(x,y)$ are two independent first integrals, where
$W_i(x,y):=\left(V_i(x,y)\right)^{2k+2}.$
\item If $\al_1=1$ then $f(x,y)=\left(1+x+y,\frac{x}{1+y}\right)$ and it preserves the
two generically transverse fibrations
$$V_1(x,y)=\frac{1-2\,x+3\,y+2\,y^2}{1+y}$$
$$V_2(x,y)=\frac{1+2\,x+3\,y+2\,y^2}{2\,(1+y)}$$
with $V_1(f(x,y))=-V_1(x,y)$ and  $V_2(f(x,y))=V_2(x,y)+1.$
Furthermore $f(x,y)$ is integrable being $W(x,y)=V_1^2(x,y)$ a first
integral.
\end{itemize}
\end{propo}
\begin{proof}

Condition $A_2 = O_0$ gives $\al_0=\ga_0=\frac{1}{\alpha_1}.$ From
Theorem \ref{theo2} we know that $f$ has two invariant fibrations.
To find them follow the procedure explained in subsection
\ref{invfib}. We consider an arbitrary cubic projective curve:
$$
\begin{array}{ll}
C[x_0:x_1:x_2]=&
{ r_0}\,{x_{{0}}}^{3}+{ r_1}\,{x_{{0}}}^{2}x_{{1}}+{ r_2}\,{x_{{0
}}}^{2}x_{{2}}+{ r_3}\,x_{{0}}{x_{{1}}}^{2}+{ r_4}\,x_{{0}}{x_{{2}
}}^{2}\\
&+{ r_5}\,x_{{0}}x_{{1}}x_{{2}}+{ r_6}\,{x_{{1}}}^{3}+{ r_7
}\,{x_{{1}}}^{2}x_{{2}}+{ r_8}\,x_{{1}}{x_{{2}}}^{2}+{ r_9}\,{x_{{
2}}}^{3}\\
\end{array}$$
and we force that $C$ is zero over the indeterminacy
points of $F,$ that is, $C(O_0)=C(O_1)=C(O_2)=0.$ Then
$$C(F^{-1}[x_0:x_1:x_2])=T_0\cdot T_1\cdot T_2\cdot
\bar{C}[x_0:x_1:x_2],$$ where $\{T_0,T_1,T_2\}=\mathcal{E}(F^{-1})$
and $\bar{C}[x_0:x_1:x_2]$ is as follows: \begin{footnotesize}$$\begin{array}{lll}
(r_2\alpha_1^2-2r_4\alpha_1+3r_9)x_0^3+(r_4\alpha_1^2-3r_9\alpha_1)x_0^2x_1+(r_2\alpha_1^3+r_1\alpha_
1^2-2r_4\alpha_1^2-r_5\alpha_1+
3r_9\alpha_1+r_8)x_0^2x_2+\\\alpha_1^2r_9x_0x_1^2+(r_5\alpha_1]2-2r_8\alpha_1)x_0x_1x_2+
(r_1\alpha_1^3-r_5\alpha_1^2+r_8\alpha_1)x_0x_2^2+(-r_9\alpha_1^3+r_8\alpha_1^2)x_1^2x_2+\\(r_9\alpha_1^3+r_3\alpha_1^2-r_3\alpha_1^2)x_1x_2^2.
\end{array}$$\end{footnotesize}
The curve $\bar{C}$ is a degree three algebraic curve. We now impose
that $C[x_0:x_1:x_2]=k\,\bar{C}[x_0:x_1:x_2],$ then after some
calculations we found (in affine coordinates)
$$\begin{array}{lll}
Q_1:=\sqrt{\alpha_1}-\alpha_1(\alpha_1+\sqrt{\alpha_1})x+\alpha_1(1+2\sqrt{\alpha_1})y+\alpha_1^2(1+\sqrt{\alpha_1})y^2,\\
Q_2:=-1+\alpha_1(1-\sqrt{\alpha_1})x+(\sqrt{\alpha_1}-2\alpha_1)y+\alpha_1(\sqrt{\alpha_1}-\alpha_1)y^2,\\
L:=1+\alpha_1  y.
\end{array}$$  The curves
$Q_1$ and $Q_2$ are invariant algebraic curves while $L$ is an
exceptional curve. Taking $V_1=Q_1/L$ and $V_2=Q_2/L,$ simple
computations prove that $V_1(f(x,y))=-\sqrt{\alpha_1}\,V_1(x,y),$
$V_2(f(x,y))=\sqrt{\alpha_1}\,V_2(x,y)$ and that $V_1(x,y),V_2(x,y)$
are generically transverse.

Now considering the mapping $\varphi(x,y):=(V_1(x,y),V_2(x,y)),$ we
see that it is a birational mapping and it has the property that
$(\varphi^{-1}\circ f \circ
\varphi)(x,y)=(-\sqrt{\al_1}x,\sqrt{\al_1}y).$ From this we deduce
that if $\al_1^{k+1}=1\,,\,\al_1\ne \pm 1$ then $f(x,y)$ is a
$(2k+2)-$periodic map. For $\al_1=-1,$ $f$ is a 4-periodic
map.Furthermore since $W_i(f(x,y))=W_i(x,y)$ for $i=1,2$ we get that
$W_1(x,y),W_2(x,y)$ are first integrals.

When $\al_1=1$ we see that $V_1$ or $V_2$ is a constant function and
that it is the unique value of the parameters which has this
behaviour. If we take $\sqrt{1}=1$ we get the invariant fibration
$V_1(x,y)={\frac{1-2\,x+3\,y+2\,y^2}{1+y}}$ with
$V_1(f(x,y))=-V_1(x,y).$ To find $V_2$ we consider a rational
function of type $V(x,y)=\frac{k_0+k_1x+k_2y+k_3y^2}{1+y}$ where
$k_i\in\C$ for $i\in \{0,1,2,3\}$ and imposing $V(f(x,y))=V(x,y)+1$ after
some calculations we find $V_2(x,y).$ Also in this case $f(x,y)$ is
birationally conjugated to $(-x,y+1),$ see Theorem
\ref{blancdeserti} again.
\end{proof}

To deal with the case $p=1,$ that is $F(A_2)=O_0,$ we notice that
this condition is equivalent to
$$\al_1^2(\al_0-\ga_0)+\al_0\al_1-1=0\quad , \quad
\al_1\ga_0(\ga_0-1)+\al_0\al_1-\ga_0=0\quad , \quad \ga_0\al_1-1\ne
0.$$ It is easy to see that it is true if and only if
$$\ga_0=\frac{1}{1+\al_1}\quad , \quad
\al_0=\frac{1+\al_1+\al_1^2}{\al_1(1+\al_1)^2}\quad , \quad
\al_1\notin\{0,-1\}.$$

We note that for these maps condition (\ref{cond.k}) reads as
$1+\al_1+\al_1^2+\cdots +\al_1^{k+1}=0,$ which implies that
$\alpha_1^{k+2}=1.$

\begin{propo}\label{pone}
Assume that $F(A_2)=O_0.$
 then $f(x,y)$ can be written as:
\begin{equation}\label{pu}f(x,y)=\left({\frac
{{\alpha_1}^{2}+\alpha_1+1}{\alpha_1\, \left( 1+\alpha_1 \right)
^{2}}} +\alpha_1\,x+y,{\frac
{x}{\frac{1}{1+\alpha_1}+y}}\right)\quad , \quad
\al_1\notin\{0,-1\}\end{equation}
%with $\alpha_1\ne
%0\,,\,\alpha_1\ne -1\quad , \quad \al_1\notin\{0,-1\}$
and \begin{itemize}
\item If $\alpha_1\ne 1$ and ${\alpha_1}^{2}+\alpha_1+1\ne 0,$
then the map $f(x,y)$ preserves the two generically transverse
fibrations
$$V_1(x,y)=\frac{B_0+B_1\,x+B_2\,y+B_3\,y^2}{C_0+C_1\,x+C_2\,y+C_3\,x\,y+C_4\,y^2}$$
$$V_2(x,y)=\frac{D_0+D_1\,x+D_2\,y+D_3\,y^2}{C_0+C_1\,x+C_2\,y+C_3\,x\,y+C_4\,y^2}$$
where

\begin{footnotesize}
$$
 \begin{array}{lll}
B_0={\alpha_1}^{2}+\alpha_1+1 & C_0={\alpha_1}^{2}+\alpha_1+1 \\
B_1=- \left( 1+\alpha_1 \right) ^{2} \left( \sqrt {\alpha_1}-\alpha_1-1 \right)\sqrt {\alpha_1} & C_1={\alpha_1}^{2} \left( 1+\alpha_1 \right) ^{2} &
\\
B_2=- \left( 1+\alpha_1 \right)  \left( {\alpha_1}^{\frac{3}{2}}-2\,{\alpha_1}^{2}-2\,\alpha_1-1 \right) & C_2=2\,{\alpha_1}^{3}+3\,
{\alpha_1}^{2}+2\,\alpha_1+1 \\
B_3=-\alpha_1\,\left(1+\alpha_1\right)^{2}\left(\sqrt{\alpha_1}-\alpha_1-1
\right) &C_3=\alpha_1\, \left( \alpha_1-1 \right)  \left( 1+\alpha_1
\right) ^{3} \\
\phantom{qweqweqweqwe} &C_4={\alpha_1}^{2} \left( 1+\alpha_1 \right)
^{2}
\end{array}
$$
\end{footnotesize}
and
\begin{footnotesize}
$$ \begin{array}{ll}
        D_0={\alpha_1}^{2}+\alpha_1+1,\\
        D_1=- \left( 1+\alpha_1 \right) ^{2} \left( \sqrt {\alpha_1}+\alpha_1+1 \right)\sqrt {\alpha_1},\\
 D_2=\left( 1+\alpha_1 \right)  \left( {\alpha_1}^{3/2}+2\,{\alpha_1}^{2}+2\,\alpha_1+1 \right),\\
D_3=\alpha_1\, \left( 1+\alpha_1 \right) ^{2} \left(
\sqrt{\alpha_1}+\alpha_1+1 \right).
\end{array},$$
\end{footnotesize}
with $V_1(f(x,y))=\frac{1}{\sqrt{\alpha_1}}\,V_1(x,y)$ and
$V_2(f(x,y))=-\frac{1}{\sqrt{\alpha_1}}\,V_2(x,y).$

If $1+\al_1+\al_1^2+\cdots +\al_1^{k+1}=0$ then $f(x,y)$ is a
$2\,(k+2)-$periodic map. In this case $W_i(x,y)$ for $i \in \{1,2\}$
are two independent first integrals, where
$$W_i(x,y):=V_i(x,y)\cdot V_i(f(x,y))\cdot V_i(f^2(x,y))\cdots
V_i(f^{2k+3}(x,y)).$$

%Furthermore $f(x,y)$ is birationally conjugated to the map
%$(\frac{1}{\sqrt{\alpha_1}}\,x,-\frac{1}{\sqrt{\alpha_1}}\,y).$

\item If $\alpha_1=1$ then $f(x,y)$ preserves the two generically
tranverse fibrations
$$V_1(x,y)={\frac {16\,xy+4\,x-6\,y-3}{4\,{y}^{2}+4\,x+8\,y+3}}$$
$$V_2(x,y)={\frac {12\,{y}^{2}-12\,x+12\,y+3}{4\,{y}^{2}+4\,x+8\,y+3}}$$
with $V_1(f(x,y))=V_1(x,y)+1$ and $V_2(f(x,y))=-V_2(x,y).$ Hence
$f(x,y)$ is integrable being $W(x,y)=V_2^2(x,y)$ a first integral.
\item
 If $\alpha_1^2+\alpha_1+1=0$ then $f(x,y)$ is a 6-periodic mapping. It preserves the two generically
tranverse fibrations
$$V_1(x,y)={\frac {2\,\alpha_1+2-x+ \left( 2\,\alpha_1-1 \right) y- \left(
1+\alpha_1
 \right) \,{y}^{2}}{ \left( \alpha_1+1 \right) x+y+ \left(
\alpha_1-1 \right) xy+ \left( \alpha_1+1 \right) {y}^{2}}}
$$
$$V_2(x,y)={\frac {\alpha_1\,x-\alpha_1\,y+{y}^{2}}{ \left( \alpha_1+1
 \right) x+y+ \left( \alpha_1-1 \right) xy+ \left( \alpha_1+1 \right) {y}^
{2}}}
$$
with $V_1(f(x,y))=-\alpha_1\,V_1(x,y)$ and
$V_2(f(x,y))=\alpha_1\,V_2(x,y).$ Furthermore $W_1(x,y):=V_1^6(x,y)$
and $W_2(x,y):=V_2^6(x,y)$ are two independent first integrals.
\end{itemize}
\end{propo}

\begin{proof}

From Theorem \ref{theo2} we know that when $1+\al_1+\al_1^2+\cdots
+\al_1^{k+1}=0,$ $d_n$ is a periodic sequence while when
$1+\al_1+\al_1^2+\cdots +\al_1^{k+1}\ne 0,$ $d_n$ is bounded. In any
case we have to find two generically transverse foliations.

We first search for invariant curves
$C(x,y)=C_0+C_1\,x+C_2\,y+C_3\,x\,y+C_4\,y^2.$ Then we consider a
rational function $V(x,y) = \frac{P(x,y)}{C(x,y)},$ where $P(x,y)$
is a second degree polynomial. The imposition of condition
$V(f(x,y))=k\cdot V(x,y)$ gives two invariant fibrations
$V_1(x,y),\,V_2(x,y)$ for
$k\in\{\frac{1}{\sqrt{\al_1}},-\frac{1}{\sqrt{\al_1}}\}.$ Also we
see that $V_1,V_2$ are generically transverse provided that
$\al_1\ne \pm 1\,,\,\al_1^2+\al_1+1\ne 0.$

Let $\varphi(x,y)$ be defined as
$\varphi(x,y)=\left(V_1(x,y),V_2(x,y)\right).$ Then $\varphi(x,y)$
is a birational map and $\varphi^{-1}\circ f \circ \varphi$ gives
the map
$\left(\frac{1}{\sqrt{\al_1}}\,x,-\frac{1}{\sqrt{\al_1}}\,y\right).$
Hence if condition (\ref{cond.k}) is accomplished, i. e., if
$1+\al_1+\al_1^2+\cdots +\al_1^{k+1}=0,$ then $f(x,y)$ is a
$(2k+4)-$periodic map.

Now assume that $\al_1=1.$ Substituting this value with $\sqrt{1}=1$
in the maps $V_1,V_2$ in the above paragraph we find that the first
fibration is a constant function while the second one is
$V_2(x,y)={\frac
{12\,{y}^{2}-12\,x+12\,y+3}{4\,{y}^{2}+4\,x+8\,y+3}},$ hence it
satisfies $V_2(f(x,y))=-V_2(x,y).$ To find the other fibration
$V(x,y)$ we consider a rational map with the same denominator of
$V_2(x,y)$ and a degree two polynomial in the numerator and imposing
$V(f(x,y))=V(x,y)+1$ we find the announced $V_1(x,y).$

The fibrations when $\al_1^2+\al_1+1=0$ are encountered in a similar
way.

\end{proof}

%Now we assume that $\tilde{F}^2(A_2)=O_0.$ It is easy to see that it
%is equivalent to $F^2(A_2)=O_0.$ For the simplification of
%calculations we consider $\al_1=\om^2.$ It implies that the
%coefficients have to satisfy:
%$$\begin{array}{l}
%Q_1:=w^6\,\ga_0^2-(\al_0\,w^6+\,(\al_0+1)\,w^4+(\al_0-2)\,w^2)\,\ga_0+w^4\,\al_0+\al_0-1=0,\\
%Q_2:=w^4\,\ga_0^3-(w^6+w^4+w^2)\,\ga_0^2+(\al_0\,w^6+(2\al_0+1)\,w^4-w^2)\,\ga_0-w^4\,\al_0-w^2\,\al_0+1=0\\
%\ga_0\,w^2-1\ne 0\quad,\quad w^2\,\ga_0^2-(w^2
%+1)\,\ga_0+w^2\,\al_0\ne 0.\end{array}
%$$
%
%Taking into account resultants of $Q_1$ and $Q_2$ we find that the condition $F(F(A_2))=O_0$ gives the maps which appears in $(a),(b),(c)$ of next theorem.
%
%The map $(b)$ corresponds to a map with $\omega=-1$ while $(c)$  corresponds to  the maps
%${\omega}^{2}-\omega+1=0$ that is, when $\al_1^2+\al_1+1=0.$
%We note that for the parametric family $(a)$ condition
%(\ref{cond.k}) is
%$$1-\om+\om^2-\om^3+\om^4+\cdots -\om^{2k+1}+\om^{2k+2}=0,$$
%which implies that $\om$ is a $(4\,k+6)-$root of unity, while for
%the  mapping $(b)$ and for the two mappings of $(c),$ condition $k$ is never
%satisfied.

\begin{propo}\label{ptwo}
     Assume that $F(F(A_2))=O_0.$
%Assume that $\tilde{F}^{2k}(A_1) \ne
%O_1$ and $\tilde{F}^{p}(A_2) = O_0$ for $(k,\,p) \in \{(k,\,2) : \,k \in \N\},$
Then $f(x,y)$ can be written as

\begin{footnotesize}$$f(x,y)=\left({\frac {{\omega}^{3}-{\omega}^{2}+1}{ \left(
\omega+1 \right)  \left( {\omega}^{2}-\omega+1 \right)
^{2}}}+{\omega}^{2}x+y,{\frac {\omega\, \left( {\omega}^{2}-\omega+1
\right) x}{\omega-1+
 \left( {\omega}^{3}-{\omega}^{2}+\omega \right) y}}
\right)\,,\, \omega\,(\omega+1)\,({\omega}^{2}-\omega+1)\ne 0,$$
\end{footnotesize}
and it preserves the fibration
$$V(x,y)=\frac{B_0+B_1\,x+B_2\,y+B_3\,y^2}{\left(\omega+ \left( {\omega}^{3}+1 \right) y\right)\left(\omega-1+ \left( {\omega}^{5}-{\omega}^{4}+{\omega}^{3}+{\omega}^{2}-
\omega+1 \right) x+ \left( {\omega}^{3}-{\omega}^{2}+\omega \right)
y \right)}$$ where
\begin{footnotesize}$$\begin{array}{ll}
        B_0=({\omega}^{3}-{\omega}^{2}+1)({\omega}-1),\\
        B_1=-{\omega}^{2} \left( \omega+1 \right)  \left( {\omega}^{2}-\omega+1
 \right) ^{2},\\
 B_2=\omega\, \left( {\omega}^{2}-\omega+1 \right)  \left( 2\,{\omega}^{3}-
{\omega}^{2}-\omega+1 \right),\\
B_3={\omega}^{3} \left( \omega+1 \right)  \left(
{\omega}^{2}-\omega+1
 \right) ^{2},
\end{array}$$\end{footnotesize}
with $V(f(x,y))=-\frac{1}{\omega}\,V(x,y).$ If $\omega^{4k+6}\ne 1$
for all $k\in\N$ this fibration is unique. If $\omega^m=(-1)^m$ for
some $m\in\N,$ then $f(x,y)$ is integrable being $W(x,y)=V(x,y)^m$ a
first integral.

When $\sum_{i=0}^{2k+2}(-1)^i\,\omega^i=0$ for a certain $k\in\N$
then $f(x,y)$ is a $(4k+6)-$periodic map.

\end{propo}

\begin{proof}

Now we assume that $F^2(A_2)=O_0.$  It is easy to see that it is
equivalent to $\tilde{F}^2(A_2)=O_0.$ For the simplification of
calculations we consider $\al_1=\om^2.$ It implies that the
coefficients have to satisfy:
$$\begin{array}{l}
E_1:=w^6\,\ga_0^2-(\al_0\,w^6+\,(\al_0+1)\,w^4+(\al_0-2)\,w^2)\,\ga_0+w^4\,\al_0+\al_0-1=0,\\
E_2:=w^4\,\ga_0^3-(w^6+w^4+w^2)\,\ga_0^2+(\al_0\,w^6+(2\al_0+1)\,w^4-w^2)\,\ga_0-w^4\,\al_0-w^2\,\al_0+1=0\\
\ga_0\,w^2-1\ne 0\quad,\quad w^2\,\ga_0^2-(w^2
+1)\,\ga_0+w^2\,\al_0\ne 0.\end{array}
$$

Taking into account some resultants of $E_1$ and $E_2$ we find that
the condition $F(F(A_2))=O_0$ gives the maps which appears in $(a).$
When ${\omega}^{2}-\omega+1=0,$ that is, when $\al_1^2+\al_1+1=0$ we
get the mappings $(b).$

We note that for the parametric family $(a)$ condition
(\ref{cond.k}) is
$$1-\om+\om^2-\om^3+\om^4+\cdots -\om^{2k+1}+\om^{2k+2}=0,$$
which implies that $\om$ is a $(4\,k+6)-$root of unity, while for
the two mappings $(b),$ condition $k$ never is satisfied.

Consider $f(x,y)$ that satisfies $(a).$ By looking for invariant curves
we find that $V(x,y)$ can be written as shown in statement of $(a).$ A
calculation shows that $V(f(x,y))=-\frac{1}{\omega}\,V(x,y).$

From this equality, we see that if $\om^m=(-1)^m$ then
$W(x,y):=V(x,y)^m$ is a first integral of $f(x,y).$

If $\sum_{i=0}^{2k+2}(-1)^i\,\omega^i=0$ for a certain $k\in\N$ then
we know that the sequence of degrees is periodic of period $4k+6.$
We are going to prove that, the map itself is periodic of
period $4k+6.$ Since $d_{4k+6} = d_0 = 1,$  the mapping $F^{4k+6}$
is linear, that is:
$$F^{4k+6}[x_0:x_1:x_2] = [r_0\,x_0+r_1\,x_1+r_2\,x_2:p_0\,x_0+p_1\,x_1+ p_2\,x_2:q_0\,x_0+q_1\,x_1+ q_2\,x_2],$$
for some constants $r_i,\,p_i,\,q_i \in \R.$ As $S_0$ is invariant
under the action $F^2,$  it is invariant under the action of
$F^{4k+6}$ as well. This implies that we can write
\begin{equation}\label{lnr-f4}
f^{4k+6}(x,y) = (p_0+p_1\,x+ p_2\,y, q_0+q_1\,x+ q_2\,y),
\end{equation}
for some $p_0,\,p_1,\,p_2,\,q_0,\,q_1,\,q_2 \in \N.$

We find that the following two are the fixed points of $f$ and the
third one is fixed by $f^2.$
$$\begin{array}{rl}
fix_1 &= \bigg(\frac{1}{(\om^2-\om+1)\,(\om+1)\,(\om^3+1)}, -\frac{\om}{\om^3+1}\bigg),\\
fix_2 &= \bigg(\frac{\om^3-\om^2+1}{\om\,(\om^2-\om+1)\,(\om^2-1)\,(\om^4-\om^3+\om-1)}, -\frac{\om^3-\om^2+1}{\om\,(\om^4-\om^3+\om-1)}\bigg),\\
fix_3 &= \bigg(\frac{1}{\om^6+2\,\om^3+1},
-\frac{\om}{\om^3+1}\bigg).\end{array}$$ Now these points must also
be fixed by $f^{4k+6}.$ Then by finding the images of $fix_1, fix_2$
and $fix_3$ under the action of $f^{4k+6}$ using (\ref{lnr-f4}) such
that $f^{4k+6}(fix_1) = fix_1,$ $f^{4k+6}(fix_2) =
fix_2,\,f^{4k+6}(fix_3) = fix_3.$ Also as the sequence of degrees is
periodic of period $4\,k+6$ this implies that
$(\tilde{F}_1^{*})^{4k+6}$ fixes the elements in the basis of Picard
group. This implies that $(\tilde{F}_1^{*})^{4k+6}$ also fixes $E_1$
that is the blown up fibre at $A_2.$ Then $F^{4k+6}$ fixes the base
point $A_2$ in $\pr.$ By utilizing this information and then solving
this system of four equations for the values of
$p_0,\,p_1,\,p_2,\,q_0,\,q_1,\,q_2$ we find that $(p_0,\,p_1 ,\,p_2
,\,q_0,\,q_1,\,q_2)=(0,1,0,0,0,1)$ which shows that $f^{4k+6}(x,y) =
(x,y).$

\end{proof}

Next case of zero entropy is when $p=3$ and $k=1.$ Condition $k$
implies $\al_1^2\,\ga_0+1=0,$ i. e., $\ga_0=\frac{-1}{\al_1^2}.$ It
is easy to see that the condition $\tilde{F}^{3}(A_2) = O_0$ is
equivalent to ${F}^{3}(A_2) = O_0.$ Some computations show that it
is true if and only if $\al_1^6+\alpha_1^3+1=0$ and
$\al_0=-2\,\al_1^5+\al_1^3-\al_1^2-\al_1.$

\begin{propo}\label{p3k1}
Assume that $F^3(A_2)=O_0$ and that condition $k$ is satisfied for
$k=1.$ Then $f(x,y)$ can be written as
$$f(x,y)=\left(-2\,\al_1^5+\al_1^3-\al_1^2-\al_1+\al_1\,x+y,\frac{x}{(\al_1+\al_1^4)+y}\right)\,,\,\al_1^6+\alpha_1^3+1=0,$$
and it is a $18-$periodic map. It preserves the two following
generically transverse foliations
$V_1(x,y)=\frac{H_1(x,y)}{C(x,y)^2}$ and
$V_2(x,y)=\frac{H_2(x,y)}{C(x,y)^2}$ where

$$C(x,y)=-{\alpha_{{1}}}^{4}-{\alpha_{{1}}}^{3}+{\alpha_{{1}}}^{2}-2-{\alpha_{{
1}}}^{4}x+ \left(
{\alpha_{{1}}}^{5}-2\,{\alpha_{{1}}}^{4}-{\alpha_{{1
}}}^{3}+2\,{\alpha_{{1}}}^{2}-\alpha_{{1}}-1 \right) y-{y}^{2}$$ and
\begin{footnotesize}$$\begin{array}{ll}
H_1(x,y)=A_0+A_1x+A_2y+A_3x^2+A_4xy+A_5y^2+A_6x^2y+A_7xy^2+A_8y^3+A_9x^3y+A_{10}x^2y^2+
\\\phantom{1231231231}A_{11}xy^3+A_{12}y^4+12x^3y^2+A_{13}x^2y^3+A_{14}xy^4,\end{array}$$\end{footnotesize}
\begin{footnotesize}$$\begin{array}{ll}
H_2(x,y)=B_0+B_1x+B_2y+B_3x^2+B_4xy+B_5y^2+B_6x^2y+B_7xy^2+B_8y^3+3\alpha_1^4x^3y+B_{9}x^2y^2+
\\\phantom{1231231231}B_{10}xy^3+B_{11}y^4+3x^2y^3+B_{12}xy^4,\end{array}$$\end{footnotesize}
with
\begin{footnotesize}
$$
 \begin{array}{lll}
A_0=31\,{\alpha_{{1}}}^{5}+23\,{\alpha_{{1}}}^{4}-23\,{\alpha_{{1}}}^{3}+
35\,\alpha_{{1}}+12
 & A_8=-4\,{\alpha_{{1}}}^{5}-16\,{\alpha_{{1}}}^{4}-4\,{\alpha_{{1}}}^{3}+10
\,{\alpha_{{1}}}^{2}-6\,\alpha_{{1}}-12\\
A_1=-18\,{\alpha_{{1}}}^{5}+4\,{\alpha_{{1}}}^{4}+20\,{\alpha_{{1}}}^{3}-
14\,{\alpha_{{1}}}^{2}-14\,\alpha_{{1}}+10&A_9=2\,{\alpha_{{1}}}^{4}+2\,{\alpha_{{1}}}^{3}+4\,\alpha_{{1}}+4\\
A_2=32\,{\alpha_{{1}}}^{5}+48\,{\alpha_{{1}}}^{4}-14\,{\alpha_{{1}}}^{3}-
20\,{\alpha_{{1}}}^{2}+44\,\alpha_{{1}}+36&A_{10}=12\,{\alpha_{{1}}}^{5}-12\,{\alpha_{{1}}}^{4}-12\,{\alpha_{{1}}}^{3}+
18\,{\alpha_{{1}}}^{2}+6\,\alpha_{{1}}-12\\
A_3=-3\,{\alpha_{{1}}}^{4}-2\,{\alpha_{{1}}}^{3}+2\,{\alpha_{{1}}}
^{2}-2&A_{11}=16\,{\alpha_{{1}}}^{5}+8\,{\alpha_{{1}}}^{4}-12\,{\alpha_{{1}}}^{3}+2
\,{\alpha_{{1}}}^{2}+16\,\alpha_{{1}}+6\\
A_4=-60\,{\alpha_{{1}}}^{5}-4\,{\alpha_{{1}}}^{4}+58\,{\alpha_{{1}}}^{3}-
34\,{\alpha_{{1}}}^{2}-48\,\alpha_{{1}}+20&A_{12}=16\,{\alpha_{{1}}}^{5}+8\,{\alpha_{{1}}}^{4}-12\,{\alpha_{{1}}}^{3}+2
\,{\alpha_{{1}}}^{2}+16\,\alpha_{{1}}+6\\
A_5=-6\,{\alpha_{{1}}}^{5}+2\,{\alpha_{{1}}}^{4}+8\,{\alpha_{{1}}}^{3}-3\,
{\alpha_{{1}}}^{2}-2\,\alpha_{{1}}+5&A_{13}=-16\,{\alpha_{{1}}}^{5}+2\,{\alpha_{{1}}}^{3}-14\,{\alpha_{{1}}}^{2}+4\\
A_6=8\,{\alpha_{{1}}}^{5}-18\,{\alpha_{{1}}}^{4}-16\,{\alpha_{{1}}}^{3}+16
\,{\alpha_{{1}}}^{2}-20&A_{14}=-4\,{\alpha_{{1}}}^{5}-4\,{\alpha_{{1}}}^{4}-2\,{\alpha_{{1}}}^{2}-2\,
\alpha_{{1}},\\
A_7=-24\,{\alpha_{{1}}}^{5}-12\,{\alpha_{{1}}}^{4}+12\,{\alpha_{{1}}}^{3}-
8\,{\alpha_{{1}}}^{2}-20\,\alpha_{{1}}&\phantom{123123123123123}
\end{array}
$$
\end{footnotesize}
and
\begin{footnotesize}
$$
 \begin{array}{lll}
B_0=-38\,{\alpha_{{1}}}^{5}-20\,{\alpha_{{1}}}^{4}+31\,{\alpha_{{1}}}^{3}-
7\,{\alpha_{{1}}}^{2}-40\,\alpha_{{1}}-7&B_7=20\,{\alpha_{{1}}}^{5}-30\,{\alpha_{{1}}}^{3}+16\,{\alpha_{{1}}}^{2}+
24\,\alpha_{{1}}-12
\\
B_1=
-11\,{\alpha_{{1}}}^{4}-4\,{\alpha_{{1}}}^{3}+9\,{\alpha_{{1}}}^{2}-4
\,\alpha_{{1}}-11 &
B_8=-3\,{\alpha_{{1}}}^{5}+7\,{\alpha_{{1}}}^{4}+3\,{\alpha_{{1}}}^{3}-9\,
{\alpha_{{1}}}^{2}+2\,\alpha_{{1}}+9\\
B_2=-62\,{\alpha_{{1}}}^{5}-51\,{\alpha_{{1}}}^{4}+44\,{\alpha_{{1}}}^{3}+
5\,{\alpha_{{1}}}^{2}-72\,\alpha_{{1}}-29
 & B_9=-3\,{\alpha_{{1}}}^{5}+6\,{\alpha_{{1}}}^{4}+6\,{\alpha_{{1}}}^{3}-6\,
{\alpha_{{1}}}^{2}+3\,\alpha_{{1}}+3\\
B_3=3\,{\alpha_{{1}}}^{5}-2\,{\alpha_{{1}}}^{3}+3\,{\alpha_{{1}}}^{2}+3\,
\alpha_{{1}}-1
 & B_{10}=3\,{\alpha_{{1}}}^{5}-9\,{\alpha_{{1}}}^{4}+6\,{\alpha_{{1}}}^{2}-6\,
\alpha_{{1}}-3\\
B_4=31\,{\alpha_{{1}}}^{5}-15\,{\alpha_{{1}}}^{4}-33\,{\alpha_{{1}}}^{3}+
32\,{\alpha_{{1}}}^{2}+21\,\alpha_{{1}}-24
 &
B_{11}=3\,{\alpha_{{1}}}^{5}+{\alpha_{{1}}}^{4}-3\,{\alpha_{{1}}}^{3}+2\,
\alpha_{{1}}\\
B_5=-25\,{\alpha_{{1}}}^{5}-20\,{\alpha_{{1}}}^{4}+19\,{\alpha_{{1}}}^{3}+
{\alpha_{{1}}}^{2}-28\,\alpha_{{1}}-10
 & B_{12}=-3\,{\alpha_{{1}}}^{5}-3\,{\alpha_{{1}}}^{2}.\\
 B_6=6\,{\alpha_{{1}}}^{5}+3\,{\alpha_{{1}}}^{4}-3\,{\alpha_{{1}}}^{3}+3\,{
\alpha_{{1}}}^{2}+6\,\alpha_{{1}}+6&\phantom{123123123123123}
\end{array}
$$
\end{footnotesize}
They satisfy $V_1(f(x,y))=\alpha_1^3V_1(x,y)$ and
$V_2(f(x,y))=\alpha_1^2V_2(x,y).$ Hence, $W_1(x,y)=V_1(x,y)^6$ and
$W_2=V_2(x,y)^9$ are two generically transverse first integrals of
$f(x,y).$
\end{propo}

\begin{proof}

To find the foliations we began looking for degree 3 invariant
curves. We only found $\bar{C}[x_0:x_1:x_2]=x_0\,C^h[x_0:x_1:x_2]$
where $C^h[x_0:x_1:x_2]$ is the homogeneous polynomial of degree two
with $C^h[1:x_1:x_2]=C(x_1,x_2).$ Then we were looking for degree
six invariant curves, with the condition that they passes trough the
three indeterminacy points $O_1,O_2$ and $O_3$ with multiplicity
two.Consequently, its image has also degree six. Forcing that this
image coincides with the curve itself we found some of them. For
instance, the two numerators of $V_1(x,y)$ and $V_2(x,y).$ A
computation gives that $V_1(f(x,y))=\alpha_1^3V_1(x,y),$
$V_2(f(x,y))=\alpha_1^2V_2(x,y)$ and that they are generically
transverse. Clearly $W_1(x,y)$ and $W_2(x,y)$ are first integrals of
$f(x,y)$ because $\alpha_1^{18}=1.$

From Theorem \ref{theo2} we know that the sequence of degrees is
periodic of periodic $18.$ To prove that the map is periodic we
apply the result of \cite{CGM2}, which says that if a map has two
independent first integrals, then it is a periodic map.

\end{proof}

%Now consider that $p=3$ and $k=2.$ Condition $k$ implies
%$\al_1^2\,\ga_0\,(1+\al_1)+1=0,$ i. e.,
%$\ga_0=\frac{-1}{\al_1^2\,(1+\al_1)}.$ Also now $\tilde{F}^{3}(A_2)
%= O_0$ is equivalent to ${F}^{3}(A_2) = O_0.$ Some tedious
%computations show that it is true if and only if either,
%$1+\al_1+\al_1^2+\alpha_1^3+\al_1^4=0$ with
%$\al_0=-\,(\al_1^3+2\,\al_1^2+\al_1+2)$ or $\al_1=1$ with
%$\al_0=\frac{1}{4}.$

\begin{propo}\label{p3k2}
Assume that $F^3(A_2)=O_0$ and that condition $k$ is satisfied for
$k=2.$
%Assume that $\tilde{F}^{2k}(A_1) = O_1$ and ${F}^{3}(A_2) =
%O_0$ and for $(k,\,p) = (2,\,3),$
Then either:
\begin{itemize}
\item [(a)] There exists $\alpha_1$ with
$\alpha_1^4+\alpha_1^3+\alpha_1^2+\alpha_1+1=0$ such that $f(x,y)$
is of the form
\begin{equation}\label{ptreskdosI}f(x,y)=\left(-(\alpha_1^3+2\,\alpha_1^2+\alpha_1+2)+\alpha_1\,x+y,\frac{x}{-(1+\alpha_1^2+\alpha_1^3)+y}
 \right).\end{equation}
 That map $f(x,y)$ preserves the elliptic fibration
 $V(x,y)=\frac{L(x,y)\cdot P(x,y)\cdot Q(x,y)}{R(x,y)^2}$ where
\begin{footnotesize}
$$\begin{array}{ll}L(x,y)=\left(-\al_1^3-2\al_1^2-2\al_1-2+(\al_1^2+\al_1)x+y\right)\\
P(x,y)= \left(yx+(-\al_1^2-1)x+\al_1^2y+\al_1^3+\al_1\right)\\
Q(x,y)=\left({\alpha_{{1}}}^{3}{y}^{2}+
\left(-{\alpha_{{1}}}^{3}-{\alpha_{{1}}}^{2}-\alpha_{{1}}-1 \right)
xy+ \left( -{\alpha_{{1}}}^{3}+{\alpha_{{1}} }^{2} \right)
y+\alpha_{{1}} \right)\\ R(x,y)=\left( {y}^{2}- \left(
3\,{\alpha_{{1}}}^{3}+3\,{\alpha_{{1}}}^{2}+2 \,\alpha_{{1}}+2
\right) y-x{\alpha_{{1}}}^{2}+{\alpha_{{1}}}^{3}-{
\alpha_{{1}}}^{2}+1 \right)\end{array}$$
\end{footnotesize}
with $V(f(x,y))=\al_1^2\,V(x,y)$ and this fibration is unique.
Furthermore $f$ is integrable being $W(x,y)=V(x,y)^5$ a first
integral of $f.$
\item [(b)]The map $f(x,y)$ is:
\begin{equation}\label{ptreskdosII} f(x,y)=\left(\frac{1}{4}+x+y,\frac{x}{-\frac{1}{2}+y}
\right).\end{equation} That map $f(x,y)$ preserves the elliptic
fibration $V(x,y)=$
\begin{footnotesize}$$ {\frac {256{x}^{3}{y}^{2}+384{x}^{2}{y}^{3}+128
x{y}^{4}+128{x}^{3}y+192{x}^{2}{y}^{2}+32x{y}^{3}-16{y}^{4}-16
{x}^{2}-8\,xy+8{y}^{2}-8x-1}{\left( -4{y}^{2}+4x+1 \right)
^{2}}}$$\end{footnotesize} with $V(f(x,y))=V(x,y)$ and this
fibration is unique. Hence $f$ is integrable.
\end{itemize}
\end{propo}

\begin{proof}
When $k=2$ condition $k$ says $\al_1^2\,\ga_0\,(1+\al_1)+1=0,$ i.
e., $\ga_0=\frac{-1}{\al_1^2\,(1+\al_1)}.$ Also here
$\tilde{F}^{3}(A_2) = O_0$ is equivalent to ${F}^{3}(A_2) = O_0.$
Some tedious computations show that it is true if and only if
either, $1+\al_1+\al_1^2+\alpha_1^3+\al_1^4=0$ with
$\al_0=-\,(\al_1^3+2\,\al_1^2+\al_1+2)$ or $\al_1=1$ with
$\al_0=\frac{1}{4}.$

For the mappings $(a)$ we find the invariant conic:
\newline $ {y}^{2}- \left(
3\,{\alpha_{{1}}}^{3}+3\,{\alpha_{{1}}}^{2}+2 \,\alpha_{{1}}+2
\right) y-{\alpha_{{1}}}^{2}\,x+{\alpha_{{1}}}^{3}-{
\alpha_{{1}}}^{2}+1$ and a degree five invariant curve, the one
given by $L(x,y)\cdot P(x,y)\cdot Q(x,y)=0.$ Taking the quotient of
them, some calculations prove that in fact
$V(f(x,y))=\al_1^2\,V(x,y).$

To prove the uniqueness of the invariant fibration we have to see
that $d_n$ is not a periodic sequence. Assume that it is, i. e.,
assume that $d_n$ is $30-$periodic. Then $F^{30}$ has degree one:
$$F^{30}\homog\,=\,[r_0\,x_0+r_1\,x_1+r_2\,x_2:p_0\,x_0+p_1\,x_1+p_2\,x_2:q_0\,x_0+q_1\,x_1+q_2\,x_2].$$
As before, since $S_0$ is invariant under $F^2$, we can write
$f^{30}$ as follows:
$$f^{30}(x,y)=\left(p_0+p_1\,x+p_2\,y,q_0+q_1\,x+q_2\,y\right).$$
Now, using that the  conic $ {y}^{2}- \left(
3\,{\alpha_{{1}}}^{3}+3\,{\alpha_{{1}}}^{2}+2 \,\alpha_{{1}}+2
\right) y-x{\alpha_{{1}}}^{2}+{\alpha_{{1}}}^{3}-{
\alpha_{{1}}}^{2}+1=0$ must be invariant under $f^{30}$ and that the
point $(-\alpha_1^3-\alpha_1^2, 1)$ (which is a fixed point for $f$)
must also be fixed for $f^{30},$ after some calculations we get that
either, $f^{30}$ is the identity or $f^{30}\circ f^{30}$ is the
identity. In any case, it would imply that $f$ is a periodic
mapping.

But we claim that the mapping $f$ itself is not periodic. If it were
the case, then $f^k(x,y)=(x,y)$ for some $k$ multiple of $30.$ We
observe that $f$ sends:
$$\{L(x,y)=0\}\longrightarrow \{P(x,y)=0\}\longrightarrow
\{Q(x,y)=0\}\longrightarrow \{L(x,y)=0\}.$$ In particular $f^3$
sends $\{Q(x,y)=0\}$ to $\{Q(x,y)=0\}.$ We see that the curve
$\{Q(x,y)=0\}$ can be parameterized by $y,$ because $\{Q(x,y)=0\}$
if and only if
$x=\varphi(y):=\frac{\al_1(\al_1^2+1+\al_1\,y)}{\al_1^2+1-y}.$ Then
$f^3(\varphi(y),y)=(\varphi(h(y)),h(y))$ where
$h(y)=\frac{u(y)}{v(y)}$ with
\begin{footnotesize}$$\begin{array}{ll}
u(y)=&-5{\alpha_{{1}}}^{3}-3{\alpha_{{1}}}^{2}-\alpha_{{1}}-6+
\left( 29
\,{\alpha_{{1}}}^{3}+10{\alpha_{{1}}}^{2}+13\alpha_{{1}}+27
 \right) y+ \left( -54\,{\alpha_{{1}}}^{3}-3\,{\alpha_{{1}}}^{2}-31\,
\alpha_{{1}}-40 \right) {y}^{2}+\\
&\left( 50{\alpha_{{1}}}^{3}-9{ \alpha_{{1}}}^{2}+32\alpha_{{1}}+20
\right) {y}^{3}+ \left( -22{
\alpha_{{1}}}^{3}+7{\alpha_{{1}}}^{2}-23\,\alpha_{{1}}-3 \right) {y}
^{4}+ \left( 2{\alpha_{{1}}}^{3}-6\,{\alpha_{{1}}}^{2}+5\alpha_{{1
}}-2 \right) {y}^{5} \end{array}$$\end{footnotesize} and

\begin{footnotesize}$$\begin{array}{ll} v(y)=&\left( 5{\alpha_{{1}}}^{3}-{\alpha_{{1}}}^{2}+4\,\alpha_{{1}}+2
 \right) y+ \left( -11\,{\alpha_{{1}}}^{3}+8{\alpha_{{1}}}^{2}-11
\alpha_{{1}} \right) {y}^{2}+ \left( 9{\alpha_{{1}}}^{3}-12\,{\alpha
_{{1}}}^{2}+12\,\alpha_{{1}}-9 \right) {y}^{3}+\\ &\left(
8\,{\alpha_{{1} }}^{2}-9\alpha_{{1}}+8 \right) {y}^{4}+ \left(
-3\,{\alpha_{{1}}}^{3 }-6\,{\alpha_{{1}}}^{2}-\alpha_{{1}}-5 \right)
{y}^{5}.\end{array}$$\end{footnotesize} If $f$ where a periodic
mapping, $h$ also would be periodic. But $h$ has the fixed point
$\bar{y}=1+\al_1^2$ and the derivative of $h(y)$ at this points
gives zero. And it is a contradiction because periodic maps have the
eigenvalues of modulus one at the fixed points.

To prove $(b)$ we begin by proving that the sequence of degrees
grows quadratically. Then the prescribed  fibration will be unique.
 The characteristic polynomial associated to
$d_n$ is $ \left( x+1 \right)  \left( {x}^{2}+x+1 \right)  \left(
{x}^{4}+{x}^{3 }+{x}^{2}+x+1 \right)  \left( x-1 \right) ^{4} $
which implies that either, $d_n$ grows quadratically or it is
periodic. It only depends on the initial conditions, that is on the
values of $d_n$ for $n=1,2,\ldots,11.$ For that mapping we have been
able to calculate these numbers: $2,3,5,8,12,16,22,28,35,43,52$
which implies that
$$d_n={\frac {97}{72}}+{\frac {5\,{n}^{2}}{12}}-\frac{1}{8}\, \left( -1 \right) ^{n}
-\frac{1}{9}\, \left(\frac{ -1+\sqrt {3}\,I}{2} \right)
^{n}-\frac{1}{9}\, \left(\frac{ -1-\sqrt {3}\,I}{2}  \right) ^{n},$$
that is, $d_n$ grows quadratically.

To find $V(x,y)$ we searched for invariant curves and we found one
of degree two: $-4{y}^{2}+4x+1$ and one of degree five, the
numerator of $V(x,y).$ Taking the quotient of them, we verified that
it satisfies $V(f(x,y))=V(x,y).$

\end{proof}

The last class with zero entropy is when $p=4$ with $k=1.$ The
condition $k=1$ says that $\ga_0=\frac{-1}{a1^2}.$ From the proof
and notations of Theorem \ref{theo2} we know that:
$$\begin{array}{ll}
S_0\longrightarrow E_0\longrightarrow S_0=T_2,\\
S_1\longrightarrow G_0\longrightarrow G_1\longrightarrow G_2
\longrightarrow T_1,\\
S_2\longrightarrow E_1\longrightarrow E_2\longrightarrow
E_3\longrightarrow E_4\longrightarrow E_5\longrightarrow T_0.
\end{array}$$
Hence, if $A_2\in S_1,$ i.e., $\al_1=-1,$ then it could happen that
$\tilde{F}^4(A_2)=O_0.$ Following the orbit of $A_2$ we get:
$\tilde{F}(A_2)=[1:0]_{G_0}\,,\,
\tilde{F}^2[1:0]_{G_0}=[1:-\al_0]_{G_2}$ and
$\tilde{F}[1:-\al_0]_{G_2}=[1:0:1]=O_0.$ Hence we see that
$\tilde{F}^4(A_2)=O_0$ for all values of $\al_0,$ provided that
$\al_1=-1$ and $\ga_0=\frac{-1}{a1^2}=-1.$

%Now assume that there are not collisions of orbits, that is,
%$\al_1+1\ne 0.$ Looking at the the expression of $F^4(A_2)$ and
%after tedious computations we get that $F^4(A_2)=O_0$ if and only if
%$f(x,y)$ is one of the following $(b)\,(c)$ or $(d)$ of next
%theorem.

\begin{propo}\label{p4k1}
Assume that $\tilde{F}^4(A_2)=O_0$ and that condition $k$ is
satisfied for $k=1,$ where $\tilde{F}$ is the mapping induced by $F$
after blowing up the point $[0:1:0].$

%Assume that $\tilde{F}^{2k}(A_1) = O_1$ and $\tilde{F}^{p}(A_2) =
%O_0$ for $(k,\,p) = (1,\,4),$
%where $\tilde{F}$ is the induced map after blowing up the point $[0:1:0].$
Then either:
\begin{itemize}
\item [(a)] The map $f(x,y)$ can be written as \begin{equation}\label{pquatrekuI}f(x,y)\,=\,\left(\alpha_0-x+y,\frac{x}{y-1}\right)\end{equation}
and it preserves the unique elliptic fibration
$$V(x,y)\,=\,\frac {\alpha_0\,xy-{x}^{2}y+x{y}^{2}}{y-1}$$
with $V(f(x,y))=V(x,y)$. Hence $f$ is integrable.

\item [(b)] The map $f(x,y)$ can be written as\begin{equation}\label{pquatrekuII}f(x,y)\,=\,\left(x+y,\frac{x}{y-1}\right)\end{equation}
and it preserves the unique elliptic fibration
$$V(x,y)\,=\,{\frac {-2\,{y}^{2}+2\,x+y+1}{xy \left( x+y \right) }
}
$$
with $V(f(x,y))=-V(x,y).$ Furthermore $f$ is integrable, being
$W(x,y)=V(x,y)^2$ a first integral of $f.$

\item[(c)] The map $f(x,y)$ can be written as\begin{equation}\label{pquatrekuIII}f(x,y)\,=\,\left(\alpha_1\,x+y,
\frac{x}{y+1}\right)\,\,\text{with}\,\,\alpha_1^2+1=0\end{equation}
and it preserves the unique elliptic fibration
$$V(x,y)\,=\,-{\frac {xy \left( \alpha_1\,y-x \right) }{ \left( \alpha_1\,y+\alpha_1-2\,x
-y-1 \right)  \left( -1+\alpha_1-2\,y \right) }}$$ with
$V(f(x,y))=\alpha_1\,V(x,y).$ Furthermore $f$ is integrable, being
$W(x,y)=V(x,y)^4$ a first integral of $f.$

\item [(d)] The map $f(x,y)$ can be written as\begin{equation}\label{pquatrekuIV}f(x,y)\,=\,\left(1-\alpha_1^3+\alpha_1\,x+y,
\frac{x}{\alpha_1^2+y}\right)\,\,\text{with}\,\,\alpha_1^4+1=0\end{equation}
and it preserves the unic elliptic fibration
$$V(x,y)\,=\,-{\frac {Q_1(x,y)\,Q_2(x,y)\,Q_3(x,y)}{\left( {\alpha_1}^{3}-1+ \left( {\alpha_1}^{2}+1 \right) x+ \left( {
\alpha_1}^{2}+\alpha_1 \right) y \right) ^{2} \left( {\alpha_1}^{3}-1+{
\alpha_1}^{2} \left( {\alpha_1}^{2}+1 \right) y \right) ^{2} }}$$ where
\begin{footnotesize}$$ \begin{array}{ll}
        Q_1(x,y)=\alpha_1^2+2\,\alpha_1+1+(2\,\alpha_1^2+\alpha_1+1)\,y+\alpha_1^3\,y^2-x\,y,\\
        Q_2(x,y)=2\,{\alpha_1}^{3}+{\alpha_1}^{2}-1+ \left( -2\,{\alpha_1}^{3}+\alpha_1+1
 \right) x+ \left( {\alpha_1}^{2}+2\,\alpha_1+1 \right) y+{\alpha_1}^{3}{x}^
{2}+{\alpha_1}^{2}xy,\\
        Q_3(x,y)=-({\alpha_1}^{2}+2\,\alpha_1+1)+ \left( 2\,{\alpha_1}^{3}+{\alpha_1}^{2}-1
 \right) y-{\alpha_1}^{2}xy.
\end{array},$$\end{footnotesize}

with $V(f(x,y))=\alpha_1^2\,V(x,y).$ Furthermore $f$ is integrable,
being $W(x,y)=V(x,y)^4$ a first integral of $f.$
\end{itemize}
\end{propo}

\begin{proof}

The mapping $(a)$ corresponds to the case $\al_1=-1$ and
$\gamma_0=-1$  when there is collisions of orbits. Looking at the
expression of $F^4(A_2)$ and after tedious computations we get that
$F^4(A_2)=O_0$ if and only if $f(x,y)$ is one of $(b),(c)$ or $(d).$

To see the uniqueness of the fibrations we have to prove that $d_n$
grows quadratically. The characteristic polynomial associated to
$d_n$ is $(x-1)^4(x+1)^2(x^2+1)(x^2+x+1)$ which implies that either,
$d_n$ grows quadratically or it is periodic. It only depends on the
initial conditions, that is on the values of $d_n$ for
$n=1,2,\ldots,10.$ For each one of the mappings which appear in the
statement, we have been able to calculate these numbers. In the four
cases they give $2,3,5,7,11,15,20,25,32,39,$ which implies that
$$d_n=\frac{23}{16}+\frac{3}{8}\,n^2-\frac{3}{16}\,(-1)^n-\frac{1}{8}\,I^n-\frac{1}{8}\,(-I)^n.$$

In order to prove $(a)$ we find the family of invariant curves $
\lambda\,(\alpha_0\,xy-{x}^{2}y+x{y}^{2})+\mu\,(y-1)=0.$ Then taking
$V=\frac{P}{Q}$ with  and $P=\alpha_0\,xy-{x}^{2}y+x{y}^{2}$ and
$Q=y-1$ we have that $V(f(x,y))=V(x,y).$

To prove $(b)$ we easily see that
$$\{x=0\}\longrightarrow \{y=0\}\longrightarrow \{x +y =0\}\longrightarrow
\{x=0\}$$ and hence $x\,y\,(y+x)$ is an invariant cubic. Then taking
$V$ as the quotient of a conic and the invariant cubic and imposing
$V(f(x,y))=k\,V(x,y)$ we found that the conic can be taken as
$-2\,{y}^{2}+2\,x+y+1$ and $k=-1.$

To prove $(c)$ we find that the straight line $
\alpha_1\,y+\alpha_1-2\,x -y-1=0$ is sent to the straight line $
-1+\alpha_1-2\,y=0$ and viceversa, which implies that their product
is an invariant curve of degree two. Also it can be seen that
$$\{x=0\}\longrightarrow \{y=0\}\longrightarrow \{\al_1\,y-x=0\}\longrightarrow
\{x=0\}$$ and hence $x\,y\,(\al_1\,y-x)$ is an invariant cubic.
Taking $V$ as the quotient of this invariant curves we get that
$V(f(x,y))=\alpha_1\,V(x,y)$ and the result follows

To see $(d)$ we began searching invariant curves of degree three and
we found (in projective coordinates) $$C[x_0:x_1:x_2]=x_0\cdot\left(
({\alpha}^{3}-1)\,x_0+ \left( {\alpha}^{2}+1 \right) x_1+ \left( {
\alpha}^{2}+\alpha \right) x_2 \right)\cdot\left(
({\alpha}^{3}-1)\,x_0+{ \alpha}^{2} \left( {\alpha}^{2}+1 \right)
x_2 \right)$$

Now we take a conic that passes through two indeterminacy points of
$F$, and we impose that its image (a conic again) also passes
through two indeterminacy points of $F.$ This gives a third conic which
we impose to be equal to the first one. With this we find
$Q_1,Q_2,Q_3.$ Then $Q_1\cdot Q_2\cdot Q_3$ is an invariant curve of
degree six. Taking $V$ as the quotient of $Q_1\cdot Q_2\cdot Q_3$
over $C[1:x:y]^2$ we get that $V(f(x,y))=\alpha_1^2\,V(x,y).$

\end{proof}

We can now state the main theorem of this section:
\begin{teo}\label{section1}
Assume that
\begin{equation}\label{eq1}
f(x,y) = \left( {\alpha _0} + {\alpha _1}x + {\alpha
    _2}y,\frac{{\beta _0} + {\beta _1}x + {\beta _2}y}{{\gamma _0} +
    {\gamma _2}y} \right)\,,\,\alpha_1\ne
0\,,\,\beta_1\ne 0\,,\,\gamma_2\ne 0\,,\,\al_2\ne 0.
\end{equation}
Then it has zero entropy if and only if after an affine change of
coordinates it can be written as one of the mappings which appear in
the statements of Propositions $5,6,7,8,9,10.$ Each one of them has
the invariant fibrations which are stated in the above propositions.
\end{teo}

\section{The subfamily $\alpha_2=0.$}
By conjugating $f(x,y)$ via $h(x,y)=\left(\be_1 \ga_2x-\frac
{\be_0\ga_2-\be_2\ga_0}{\be_1\ga_2},\be_1 y-\frac
{\ga_0}{\ga_2}\right)$ and renaiming the parameters, we can consider
$$f(x,y)=\left(\al_0+\al_1 x,\frac{x+\be_2 y}{y}\right)\,\text{with}\,\, \al_1\ne 0.$$
We consider the induced map in the projective plane : $F:\pr
\rightarrow \pr$ given by
\begin{equation}\label{F1}
F[x_0:x_1:x_2]=[x_0\,x_2:(\al_0\,x_0+\al_1\,x_1)\,x_2:x_0\,(x_1+\be_2\,x_2)].
\end{equation}
The indeterminacy sets of $F$ and $F^{-1}$ are $\mathcal{I}(F) =
\left\{ {{O_0},{O_1},{O_2}} \right\},$ where
\begin{equation*}
{O_0} = \left[ 1:0:0 \right],\quad {O_1} = \left[ { 0: 0: 1}
\right], \quad {O_2} = \left[0: 1: 0\right],
\end{equation*}
and $\mathcal{I}(F^{-1}) = \left\{ {{A_0},{A_1},{A_2}} \right\},$
where
\begin{equation*}
{A_0} = \left[ 0 : 1 : 0\right],\quad {A_1} = \left[ { 0: 0 : 1}
\right], \quad {A_2} = \left[1 : \al_0 : \be_2 \right].
\end{equation*}
Furthermore the exceptional curves of $F$ and $F^{-1}$ are the
following:
\begin{equation*}
{S_0} = \left\{ x_0= 0\right\},\quad {S_1} = \left\{x_2 = 0
\right\}, \quad {S_2} = \left \{x_1 = 0\right\},
\end{equation*}
\begin{equation*}
{T_0} = \left\{ \al_0 x_0 - x_1= 0\right\},\quad {T_1} =
\left\{\be_2 x_0-x_2 = 0 \right\}, \quad {T_2} = \left \{x_0 =
0\right\}.
\end{equation*}

\begin{teo}\label{theo1}
Let $f(x,y)$ be a  map of type (\ref{eq1}) with $\gamma_2\ne 0,
\alpha_1\ne 0, \beta_1\ne 0$ and suppose that $\alpha_2 = 0$. If
$f^p(\alpha_0,\beta_2)=(0,0)$ for some $p \in \N$ then the
characteristic polynomial associated with $f$ is given by
\begin{equation*}
\mathcal{X}_p = (x^{p+1}+1)(x-1)^2(x+1),
 \end{equation*}
and the sequence of degrees of $f$ is periodic with period $2p+2.$
If no such $p$ exists then the characteristic polynomial associated
with $f$ is
$$\mathcal{X}\,=\,(x-1)^2\,(x+1),$$ and the sequence of degrees $d_n$ grows
linearly.
\end{teo}

\begin{proof}
Observe that $S_0 \twoheadrightarrow A_0 = O_2$ and
$S_1\twoheadrightarrow A_1 = O_1.$ Hence we blow up the points
$A_0,\,A_1$  getting the exceptional fibres $E_0,\,E_1.$ Let $X$ be
the new space and let $\tilde{F}: X \to X$ be the corresponding map
on $X$. Then the map $\tilde{F}$ sends the curve $S_0\to E_0 \to
S_0$ and $S_1 \to E_1\to T_1.$ We observe that no new indeterminacy
points are created therefore $\mathcal{{I}}(\tilde{F}) = \left\{
{{O_0}} \right\}$ and $\mathcal{{E}}(\tilde{F}) = \left\{ {{S_2}}
\right\}.$

Assume that there exists $p\in\N$ such that $\tilde{F}^p(A_2)=O_0.$
Then we blow up $A_2,\tilde{F}(A_2),$
$\tilde{F^2}(A_2),\ldots,\tilde{F^p}(A_2)=O_0$ getting the
exceptional fibres which we call $E_2,E_3,\ldots,E_{p+2}.$ Set
$\tilde{F_1}:X_1 \to X_1$ the extended map. Performing the blow up
at $O_0,$ since $T_0$ is sent to $O_0$ via $F^{-1},$ we have that
$\tilde{F}_1^{-1}:T_0\to E_{p+2}.$ Then
 $S_2 \to E_2 \to E_3 \to \cdots\to E_{p+1}\to E_{p+2}
\to T_0.$ Hence $\tilde{F}_1:X_1\to X_1$ is an AS map and also an
automorphism. Taking into account that $A_2=[1:\al_0:\be_2]$ and
$O_0=[1:0:0]$ belong to the affine plane, it is clear that condition
$\tilde{F}^p(A_2)=O_0$ reads as $f^p(\alpha_0,\beta_2)=(0,0).$

Now we have two closed lists as follows
$$\mathcal{L}_{c_1} = \{\mathcal{O}_0 = \{A_0 = O_2\},\quad \mathcal{O}_2 = \{A_2,\,\tilde{F}(A_2)\,,\ldots,\tilde{F}^{p}(A_2)=O_0\}\},$$
$$\mathcal{L}_{c_2} = \{\mathcal{O}_1 = \{A_1 = O_1\}\}.$$
Then by using Theorem \ref{th_BK} we find that the characteristic
polynomial associated to $F$ is $$\mathcal{X} =
(x^{p+1}+1)(x-1)^2(x+1).$$ If $p$ is even then $x^{p+1}+1$ has the
factor $x+1$ and $\mathcal{X} =
(x-1)^2\,(x+1)^2\,(x^p-x^{p-1}+\cdots -x+1).$ Hence the sequence of
degrees is $d_{n} = c_0+c_1\, n+
c_2\,(-1)^n+c_3\,n\,(-1)^n+c_4\,\lambda_1^n + c_5\, \lambda_2^n
+...+c_{p+3}\,\lambda_{p}^n,$ where $c_i$ are constants and
$\lambda_1,\,\lambda_2,...,\lambda_{p}$ are the roots of polynomial
$x^p-x^{p-1} +\cdots -x+1.$ By looking at $d_n$ we see that $f$ does
not grow quadratically or exponentially. As our map $\tilde{F}_1$ is
an automorphism then by using the results from Diller and Favre in
\cite{DF} we see that also cannot have linear growth. Therefore we
must have $c_1=c_3=0.$ Hence the sequence of degrees must be
periodic. This implies that $d_{2p+2+n} = d_{n}$ i.e. the sequence
of degrees is periodic with period $2p+2.$ If $p$ is odd then $d_n$
is also periodic of period $2p+2.$

If $\tilde{F}^{p}(A_2)\ne O_0$ for all $p\in\N,$ then we have two
lists which are open and closed as follows:
$$\mathcal{L}_{o} = \{ \mathcal{O}_{0} = \{A_0 = O_2\}\}\quad,\quad \mathcal{L}_{c} = \{ \mathcal{O}_{1} = \{A_1 = O_1\}\}.$$ Then $\delta(F)$ is determined by the polynomial $(x-1)^2(x+1),$ and $\delta(f)=1.$
The sequence of degrees is
$d_n=\frac{5}{4}+\frac{1}{2}\,n-\frac{1}{4}\,(-1)^n.$
\end{proof}

\begin{teo}\label{section2}
    Let $f(x,y)=\left(\al_0+\al_1 x,\frac{x+\be_2 y}{y}\right)$ with $\alpha_1\ne
0$ and set $h(x)=\al_0+\al_1 x.$ Then the following hold:
    \begin{enumerate}
    \item If $f^p(\alpha_0,\beta_2)\ne (0,0)$ for all $p\in\N$ then
    $f$ preserves the fibration $V_1(x,y) = x$ with $V_1(f(x,y)) = \al_0+\al_1 V_1(x,y),$ and this fibration is unique.
    If $\alpha_1^n=1$ for some $n>1\,,\,\alpha_1\ne 1,$ the map is
    integrable being $$W(x,y)=x\cdot h(x)\cdot h(h(x))\cdots  h^{n-1}(x)$$ a first integral of $f.$
    Also when $\alpha_1=1$ and $\alpha_0=0,$ $f$ is integrable.
    \item If $f^p(\alpha_0,\beta_2)=(0,0)$ for some $p\ge 1,$ then
    $f$ is a $(2p+2)-$periodic map. These maps have $W(x,y)=x\cdot h(x)\cdot h(h(x))\cdots
    h^{2p+1}(x)$ as a first integral.
    \item If $(\alpha_0,\beta_2)=(0,0),$ then $f(x,y)=(\alpha_1 x,\frac{x}{y})$
    and it preserves the two generically transverse fibrations
    $$V_1(x,y)=\sqrt{\alpha_1}y+\frac{x}{y}\quad,\quad V_2(x,y)=-\sqrt{\alpha_1}y+\frac{x}{y}$$
    with $V_1(f(x,y)) = \sqrt{\alpha_1}\,V_1(x,y)$ and
    $V_2(f(x,y))=-\sqrt{\alpha_1}\,V_2(x,y).$
    %Furthermore $f(x,y)$ is birationally conjugated to the map
    %$(\sqrt{\alpha_1}x,-\sqrt{\alpha_1}y).$
When $\alpha_1^n=1$ for some $n$ then $f$ is $2n-$periodic and
$W_1(x,y)=V_1(x,y)^{2n},W_2(x,y)=V_2(x,y)^{2n}$ are two independent
first integrals.
    \end{enumerate}

\end{teo}

\begin{remark} We notice that when $p=0$, that is, $\alpha_0=0=\beta_2,$
then $\varphi(x,y):=$\linebreak $=(V_1(x,y),V_2(x,y))$ is a
birational map. It turns out that using $\varphi(x,y)$ as a
conjugation we get the map $(\sqrt{\alpha_1}x,-\sqrt{\alpha_1}y).$
These result on linearizations was already pointed out on the work
of Blanc and Deserti, see\cite{BD}. Furthermore, the sequence of
degrees is $d_n=2,1,2,1,2,1,...$ a two-periodic sequence, and
avoiding the case $\alpha_1^n=1$ for some $n,$ the map itself in not
more periodic.

For $p\ge 1$ the map is periodic and hence it has two independent
first integrals. There is a method to find them (see \cite{CGM2}).
For instance, when $\al_1=-1$ and $\al_0=-\be_2^2$ (case $p=1$)
i.e., $f(x,y)=\left(-x-\be_2^2,\frac{x+\be_2 y}{y}\right),$ we have
that
$$H(x,y)=y+\frac{x+\be_2 y}{y}+\frac{x\,(\be_2-y)}{x+\be_2 y}+\frac{x+\be_2^2}{\be_2-y}$$
is a first integral of $f$ and $W(x,y)\,,\,H(x,y)$ are generically
transverse.
\end{remark}

\begin{proof}
If $f^p(\alpha_0,\beta_2)\ne (0,0)$ for all $p\in\N$ then from the
above theorem we know that $d_n$ grows linearly, and hence we know
that $f(x,y)$ has a unique invariant fibration. Clearly $V_1(x,y)=x$
is an invariant fibration and when $\al_1=1$ and $\al_0=0,$
$V_1(x,y)$ is a first integral. When $\al_1^n=1$ the function $h(x)$
is periodic of period $n$ and hence $W(x,y)$ is a first integral of
$f(x,y).$

Now assume that $f^p(\alpha_0,\beta_2)=(0,0)$ for a certain
$p\in\N.$ From Theorem (\ref{theo1}) we know that the sequence of
degrees $d_n$ is $2p+2$ periodic. We are going to see that $f(x,y)$
itself is a periodic map. Since the map $F^{2p+2}$ is linear, we can
consider that for some constants $r_i,\,p_i,\,q_i \in \R$ the map
$F^{2p+2}$ can be written in the following form:
$$F^{2p+2}[x_0:x_1:x_2] = [r_0\,x_0+r_1\,x_1+r_2\,x_2:p_0\,x_0+p_1\,x_1+ p_2\,x_2:q_0\,x_0+q_1\,x_1+ q_2\,x_2].$$

We know that $S_0$ is invariant under the action $F^2$ therefore it
is invariant under the action of $F^{2p+2}$ as well. This implies
that

$$F^{2p+2}[0:x_1:x_2] = [0:x_1:x_2],$$ which further implies that $r_1\,x_1+r_2\,x_2 = 0$ for all complex numbers $x_1,\,x_2.$
This is only possible if $r_1 = r_2 = 0.$ Then we can write

$$F^{2p+2}[x_0:x_1:x_2] = \bigg[x_0:\frac{p_0}{r_0}\,x_0+\frac{p_1}{r_0}\,x_1+ \frac{p_2}{r_0}\,x_2:\frac{q_0}{r_0}\,x_0+\frac{q_1}{r_0}\,x_1+\frac{q_2}{r_0}\,x_2\bigg],$$
which in the affine plane by taking $x_0 = 1$ and rewriting the
parameters, as new parameters, the function $F^{2p+2}$ can be
written as following:

\begin{equation}\label{lnr-fff}
f^{2p+2}(x,y) = (p_0+p_1\,x+ p_2\,y, q_0+q_1\,x+ q_2\,y),
\end{equation}
for any $p_0,\,p_1,\,p_2,\,q_0,\,q_1,\,q_2 \in \R.$ We find that the
following two points are fixed for $f(x,y):$
$$
(X,\pm
Y)=\bigg(\frac{\al_0}{1-\al_1},\frac{\be_2(1-\al_1)\pm\sqrt{(1-\al_1)^2\,\be_2^2+4\al_0(1-\al_1)}}{2(1-\al_1)}
\bigg),
$$

As these points are fixed by $f$ so they are also fixed points of
$f^{2p+2}.$ Then finding their images under the action of $f^{2p+2}$
using (\ref{lnr-fff}) we get a system equations such that
$f^{2p+2}(X,Y)[1]=X,\,f^{2p+2}(X,Y)[2]=Y,\,\,f^{2p+2}(X,-Y)[2]=-Y.$
Also as the sequence of degrees is periodic of period $2p+2$ this
implies that $(\tilde{F}_1^{*})^{2p+2}$ fixes the elements in the
basis of Picard group. This implies that $(\tilde{F}_1^{*})^{2p+2}$
also fixes $E_2$ that is the blown up fiber at $A_2.$ Then
$F^{2p+2}$ fixes the base point $A_2$ in $\pr.$ By utilizing this
information and then solving the system of equations for the values
of $p_0,\,p_1,\,p_2,\,q_0,\,q_1,\,q_2$ we find that $(p_0,\,p_1
,\,p_2 ,\,q_0,\,q_1,\,q_2)=(0,1,0,0,0,1)$ which implies that
$f^{2p+2}(x,y) = (x,y).$

Finally, $p = 0$ that is when $A_2 = O_0,$ by iterating the function
$f(x,y)= \left(\al_1\,x, \frac{x}{y} \right)$ we find that
$f^{2n}(x,y)= (\al_1^{2n}\,x, \al_1^{n}\,y)$ and $f^{2n+1}(x,y)=
(\al_1^{2n+1}\,x, \al_1^{n}\,\frac{x}{y}).$ Now observe that for
$\al_1^n = 1$ we have $f^{2n}(x,y)= (x,y)$ and $f^{2n+1}(x,y)=
\left(x, \frac{x}{y} \right).$ Therefore $f$ is $2n-$periodic. Now
for $\al_1^n \neq 1$ through simple calculations we find that $f$
preserves the announced fibrations $V_1(x,y)$ and $V_2(x,y).$
\end{proof}

\end{document}